\documentclass[sn-mathphys-num]{sn-jnl}

\usepackage{graphicx}%
\usepackage{multirow}%
\usepackage{amsmath,amssymb,amsfonts}%
\usepackage{amsthm}%
\usepackage{mathrsfs}%
\usepackage[title]{appendix}%
\usepackage{xcolor}%
\usepackage{textcomp}%
\usepackage{manyfoot}%
\usepackage{booktabs}%
\usepackage{algorithm}%
\usepackage{algorithmicx}%
\usepackage{algpseudocode}%
\usepackage{listings}%
\newtheorem{theorem}{Theorem}

\newtheorem{axiom}[theorem]{Axiom}

\newtheorem{conjecture}[theorem]{Conjecture}

\newtheorem{example}[theorem]{Example}
\newtheorem{exercise}[theorem]{Exercise}

\newtheorem{remark}[theorem]{Remark}

\makeatletter

\ifx\pdfoutput\relax\let\pdfoutput=\undefined\fi
\newcount\msipdfoutput
\ifx\pdfoutput\undefined
\else
 \ifcase\pdfoutput
 \else
    \ifx\paperwidth\undefined
    \else
      \ifdim\paperheight=0pt\relax
      \else
        \pdfpageheight\paperheight
      \fi
      \ifdim\paperwidth=0pt\relax
      \else
        \pdfpagewidth\paperwidth
      \fi
    \fi
  \fi
\fi

\newcount\@hour\newcount\@minute\chardef\@x10\chardef\@xv60
\def\tcitime{
\def\@time{%
  \@minute\time\@hour\@minute\divide\@hour\@xv
  \ifnum\@hour<\@x 0\fi\the\@hour:%
  \multiply\@hour\@xv\advance\@minute-\@hour
  \ifnum\@minute<\@x 0\fi\the\@minute
  }}%

\def\x@hyperref#1#2#3{%
   \catcode`\~ = 12
   \catcode`\$ = 12
   \catcode`\_ = 12
   \catcode`\# = 12
   \catcode`\& = 12
   \catcode`\% = 12
   \y@hyperref{#1}{#2}{#3}%
}

\def\y@hyperref#1#2#3#4{%
   #2\ref{#4}#3
   \catcode`\~ = 13
   \catcode`\$ = 3
   \catcode`\_ = 8
   \catcode`\# = 6
   \catcode`\& = 4
   \catcode`\% = 14
}

\@ifundefined{hyperref}{\let\hyperref\x@hyperref}{}
\@ifundefined{msihyperref}{\let\msihyperref\x@hyperref}{}

\@ifundefined{qExtProgCall}{\def\qExtProgCall#1#2#3#4#5#6{\relax}}{}
\def\QCTOpt[#1]#2{%
  \def\QCTOptB{#1}
  \def\QCTOptA{#2}
}
\def\QCTNOpt#1{%
  \def\QCTOptA{#1}
  \let\QCTOptB\empty
}
\def\Qct{%
  \@ifnextchar[{%
    \QCTOpt}{\QCTNOpt}
}
\def\QCBOpt[#1]#2{%
  \def\QCBOptB{#1}%
  \def\QCBOptA{#2}%
}
\def\QCBNOpt#1{%
  \def\QCBOptA{#1}%
  \let\QCBOptB\empty
}
\def\Qcb{%
  \@ifnextchar[{%
    \QCBOpt}{\QCBNOpt}%
}
\def\PrepCapArgs{%
  \ifx\QCBOptA\empty
    \ifx\QCTOptA\empty
      {}%
    \else
      \ifx\QCTOptB\empty
        {\QCTOptA}%
      \else
        [\QCTOptB]{\QCTOptA}%
      \fi
    \fi
  \else
    \ifx\QCBOptA\empty
      {}%
    \else
      \ifx\QCBOptB\empty
        {\QCBOptA}%
      \else
        [\QCBOptB]{\QCBOptA}%
      \fi
    \fi
  \fi
}
\newcount\GRAPHICSTYPE
\GRAPHICSTYPE=\z@
\def\GRAPHICSPS#1{%
 \ifcase\GRAPHICSTYPE
   \special{ps: #1}%
 \or
   \special{language "PS", include "#1"}%
 \fi
}%
\def\graffile#1#2#3#4{%
    \bgroup
	   \@inlabelfalse
       \leavevmode
       \@ifundefined{bbl@deactivate}{\def~{\string~}}{\activesoff}%
        \raise -#4 \BOXTHEFRAME{%
           \hbox to #2{\raise #3\hbox to #2{\null #1\hfil}}}%
    \egroup
}%
\def\draftbox#1#2#3#4{%
 \leavevmode\raise -#4 \hbox{%
  \frame{\rlap{\protect\tiny #1}\hbox to #2%
   {\vrule height#3 width\z@ depth\z@\hfil}%
  }%
 }%
}%
\newcount\@msidraft
\@msidraft=\z@
\let\nographics=\@msidraft
\newif\ifwasdraft
\wasdraftfalse

\def\GRAPHIC#1#2#3#4#5{%
   \ifnum\@msidraft=\@ne\draftbox{#2}{#3}{#4}{#5}%
   \else\graffile{#1}{#3}{#4}{#5}%
   \fi
}
\def\addtoLaTeXparams#1{%
    \edef\LaTeXparams{\LaTeXparams #1}}%
\newif\ifBoxFrame \BoxFramefalse
\newif\ifOverFrame \OverFramefalse
\newif\ifUnderFrame \UnderFramefalse

\def\BOXTHEFRAME#1{%
   \hbox{%
      \ifBoxFrame
         \frame{#1}%
      \else
         {#1}%
      \fi
   }%
}

\def\doFRAMEparams#1{\BoxFramefalse\OverFramefalse\UnderFramefalse\readFRAMEparams#1\end}%
\def\readFRAMEparams#1{%
 \ifx#1\end%
  \let\next=\relax
  \else
  \ifx#1i\dispkind=\z@\fi
  \ifx#1d\dispkind=\@ne\fi
  \ifx#1f\dispkind=\tw@\fi
  \ifx#1t\addtoLaTeXparams{t}\fi
  \ifx#1b\addtoLaTeXparams{b}\fi
  \ifx#1p\addtoLaTeXparams{p}\fi
  \ifx#1h\addtoLaTeXparams{h}\fi
  \ifx#1X\BoxFrametrue\fi
  \ifx#1O\OverFrametrue\fi
  \ifx#1U\UnderFrametrue\fi
  \ifx#1w
    \ifnum\@msidraft=1\wasdrafttrue\else\wasdraftfalse\fi
    \@msidraft=\@ne
  \fi
  \let\next=\readFRAMEparams
  \fi
 \next
 }%
\def\IFRAME#1#2#3#4#5#6{%
      \bgroup
      \let\QCTOptA\empty
      \let\QCTOptB\empty
      \let\QCBOptA\empty
      \let\QCBOptB\empty
      #6%
      \parindent=0pt
      \leftskip=0pt
      \rightskip=0pt
      \setbox0=\hbox{\QCBOptA}%
      \@tempdima=#1\relax
      \ifOverFrame
          \typeout{This is not implemented yet}%
          \show\HELP
      \else
         \ifdim\wd0>\@tempdima
            \advance\@tempdima by \@tempdima
            \ifdim\wd0 >\@tempdima
               \setbox1 =\vbox{%
                  \unskip\hbox to \@tempdima{\hfill\GRAPHIC{#5}{#4}{#1}{#2}{#3}\hfill}%
                  \unskip\hbox to \@tempdima{\parbox[b]{\@tempdima}{\QCBOptA}}%
               }%
               \wd1=\@tempdima
            \else
               \textwidth=\wd0
               \setbox1 =\vbox{%
                 \noindent\hbox to \wd0{\hfill\GRAPHIC{#5}{#4}{#1}{#2}{#3}\hfill}\\%
                 \noindent\hbox{\QCBOptA}%
               }%
               \wd1=\wd0
            \fi
         \else
            \ifdim\wd0>0pt
              \hsize=\@tempdima
              \setbox1=\vbox{%
                \unskip\GRAPHIC{#5}{#4}{#1}{#2}{0pt}%
                \break
                \unskip\hbox to \@tempdima{\hfill \QCBOptA\hfill}%
              }%
              \wd1=\@tempdima
           \else
              \hsize=\@tempdima
              \setbox1=\vbox{%
                \unskip\GRAPHIC{#5}{#4}{#1}{#2}{0pt}%
              }%
              \wd1=\@tempdima
           \fi
         \fi
         \@tempdimb=\ht1
         \advance\@tempdimb by -#2
         \advance\@tempdimb by #3
         \leavevmode
         \raise -\@tempdimb \hbox{\box1}%
      \fi
      \egroup%
}%
\def\DFRAME#1#2#3#4#5{%
  \vspace\topsep
  \hfil\break
  \bgroup
     \leftskip\@flushglue
	 \rightskip\@flushglue
	 \parindent\z@
	 \parfillskip\z@skip
     \let\QCTOptA\empty
     \let\QCTOptB\empty
     \let\QCBOptA\empty
     \let\QCBOptB\empty
	 \vbox\bgroup
        \ifOverFrame
           #5\QCTOptA\par
        \fi
        \GRAPHIC{#4}{#3}{#1}{#2}{\z@}%
        \ifUnderFrame
           \break#5\QCBOptA
        \fi
	 \egroup
  \egroup
  \vspace\topsep
  \break
}%
\def\FFRAME#1#2#3#4#5#6#7{%
  \@ifundefined{floatstyle}
    {
     \begin{figure}[#1]%
    }
    {
	 \ifx#1h
      \begin{figure}[H]%
	 \else
      \begin{figure}[#1]%
	 \fi
	}
  \let\QCTOptA\empty
  \let\QCTOptB\empty
  \let\QCBOptA\empty
  \let\QCBOptB\empty
  \ifOverFrame
    #4
    \ifx\QCTOptA\empty
    \else
      \ifx\QCTOptB\empty
        \caption{\QCTOptA}%
      \else
        \caption[\QCTOptB]{\QCTOptA}%
      \fi
    \fi
    \ifUnderFrame\else
      \label{#5}%
    \fi
  \else
    \UnderFrametrue%
  \fi
  \begin{center}\GRAPHIC{#7}{#6}{#2}{#3}{\z@}\end{center}%
  \ifUnderFrame
    #4
    \ifx\QCBOptA\empty
      \caption{}%
    \else
      \ifx\QCBOptB\empty
        \caption{\QCBOptA}%
      \else
        \caption[\QCBOptB]{\QCBOptA}%
      \fi
    \fi
    \label{#5}%
  \fi
  \end{figure}%
 }%
\newcount\dispkind%

\def\makeactives{
  \catcode`\"=\active
  \catcode`\;=\active
  \catcode`\:=\active
  \catcode`\'=\active
  \catcode`\~=\active
}
\bgroup
   \makeactives
   \gdef\activesoff{%
      \def"{\string"}%
      \def;{\string;}%
      \def:{\string:}%
      \def'{\string'}%
      \def~{\string~}%
    }
\egroup

\def\FRAME#1#2#3#4#5#6#7#8{%
 \bgroup
 \ifnum\@msidraft=\@ne
   \wasdrafttrue
 \else
   \wasdraftfalse%
 \fi
 \def\LaTeXparams{}%
 \dispkind=\z@
 \def\LaTeXparams{}%
 \doFRAMEparams{#1}%
 \ifnum\dispkind=\z@\IFRAME{#2}{#3}{#4}{#7}{#8}{#5}\else
  \ifnum\dispkind=\@ne\DFRAME{#2}{#3}{#7}{#8}{#5}\else
   \ifnum\dispkind=\tw@
    \edef\@tempa{\noexpand\FFRAME{\LaTeXparams}}%
    \@tempa{#2}{#3}{#5}{#6}{#7}{#8}%
    \fi
   \fi
  \fi
  \ifwasdraft\@msidraft=1\else\@msidraft=0\fi{}%
  \egroup
 }%
\def\TEXUX#1{"texux"}

\def\limfunc#1{\mathop{\rm #1}}%
\def\func#1{\mathop{\rm #1}\nolimits}%
\long\def\QQQ#1#2{%
     \long\expandafter\def\csname#1\endcsname{#2}}%
\@ifundefined{QTP}{\def\QTP#1{}}{}
\@ifundefined{QEXCLUDE}{\def\QEXCLUDE#1{}}{}
\@ifundefined{Qlb}{}{}
\@ifundefined{Qlt}{}{}
\long\def\QQA#1#2{}%
\def\QTR#1#2{{\csname#1\endcsname {#2}}}%

\def\EXPAND#1[#2]#3{}%
\def\NOEXPAND#1[#2]#3{}%
\def\LaTeXparent#1{}%
\def\ChildStyles#1{}%
\def\ChildDefaults#1{}%
\def\QTagDef#1#2#3{}%

\@ifundefined{correctchoice}{}{}
\@ifundefined{HTML}{\def\HTML#1{\relax}}{}
\@ifundefined{TCIIcon}{\def\TCIIcon#1#2#3#4{\relax}}{}
\if@compatibility
\else
  \providecommand{\UNICODE}[2][]{\protect\rule{.1in}{.1in}}
  \providecommand{\U}[1]{\protect\rule{.1in}{.1in}}
  
\fi

\@ifundefined{lambdabar}{
      
   }{}

\@ifundefined{StyleEditBeginDoc}{}{}
\def\QQfnmark#1{\footnotemark}

\@ifundefined{TCIMAKEINDEX}{}{\makeindex}%
\@ifundefined{abstract}{%
 \def\abstract{%
  \if@twocolumn
   \section*{Abstract (Not appropriate in this style!)}%
   \else \small
   \begin{center}{\bf Abstract\vspace{-.5em}\vspace{\z@}}\end{center}%
   \quotation
   \fi
  }%
 }{%
 }%
\@ifundefined{endabstract}{\def\endabstract
  {\if@twocolumn\else\endquotation\fi}}{}%
\@ifundefined{maketitle}{\def\maketitle#1{}}{}%
\@ifundefined{affiliation}{\def\affiliation#1{}}{}%
\@ifundefined{proof}{}{}%
\@ifundefined{endproof}{}{}%
\@ifundefined{newfield}{\def\newfield#1#2{}}{}%
\@ifundefined{chapter}{\def\chapter#1{\par(Chapter head:)#1\par }%
 \newcount\c@chapter}{}%
\@ifundefined{part}{\def\part#1{\par(Part head:)#1\par }}{}%
\@ifundefined{section}{\def\section#1{\par(Section head:)#1\par }}{}%
\@ifundefined{subsection}{\def\subsection#1%
 {\par(Subsection head:)#1\par }}{}%
\@ifundefined{subsubsection}{\def\subsubsection#1%
 {\par(Subsubsection head:)#1\par }}{}%
\@ifundefined{paragraph}{\def\paragraph#1%
 {\par(Subsubsubsection head:)#1\par }}{}%
\@ifundefined{subparagraph}{\def\subparagraph#1%
 {\par(Subsubsubsubsection head:)#1\par }}{}%
\@ifundefined{therefore}{}{}%
\@ifundefined{backepsilon}{}{}%
\@ifundefined{yen}{}{}%
\@ifundefined{registered}{%
   \def\registered{\relax\ifmmode{}\r@gistered
                    \else$\m@th\r@gistered$\fi}%
 \def\r@gistered{^{\ooalign
  {\hfil\raise.07ex\hbox{$\scriptstyle\rm\text{R}$}\hfil\crcr
  \mathhexbox20D}}}}{}%
\@ifundefined{Eth}{}{}%
\@ifundefined{eth}{}{}%
\@ifundefined{Thorn}{}{}%
\@ifundefined{thorn}{}{}%
\@ifundefined{degree}{}{}%
\newdimen\theight
\@ifundefined{Column}{\def\Column{%
 \vadjust{\setbox\z@=\hbox{\scriptsize\quad\quad tcol}%
  \theight=\ht\z@\advance\theight by \dp\z@\advance\theight by \lineskip
  \kern -\theight \vbox to \theight{%
   \rightline{\rlap{\box\z@}}%
   \vss
   }%
  }%
 }}{}%
\@ifundefined{qed}{\def\qed{%
 \ifhmode\unskip\nobreak\fi\ifmmode\ifinner\else\hskip5\p@\fi\fi
 \hbox{\hskip5\p@\vrule width4\p@ height6\p@ depth1.5\p@\hskip\p@}%
 }}{}%
\@ifundefined{cents}{}{}%
\@ifundefined{tciLaplace}{}{}%
\@ifundefined{tciFourier}{}{}%
\@ifundefined{textcurrency}{}{}%
\@ifundefined{texteuro}{}{}%
\@ifundefined{euro}{}{}%
\@ifundefined{textfranc}{}{}%
\@ifundefined{textlira}{}{}%
\@ifundefined{textpeseta}{}{}%
\@ifundefined{miss}{\def\miss{\hbox{\vrule height2\p@ width 2\p@ depth\z@}}}{}%
\@ifundefined{vvert}{}{}
\@ifundefined{tcol}{\def\tcol#1{{\baselineskip=6\p@ \vcenter{#1}} \Column}}{}%
\@ifundefined{dB}{}{}
\@ifundefined{mB}{}{}
\@ifundefined{nB}{}{}
\@ifundefined{note}{}{}%
\def\newfmtname{LaTeX2e}
\ifx\fmtname\newfmtname
  \DeclareOldFontCommand{\rm}{\normalfont\rmfamily}{\mathrm}
  \DeclareOldFontCommand{\sf}{\normalfont\sffamily}{\mathsf}
  \DeclareOldFontCommand{\tt}{\normalfont\ttfamily}{\mathtt}
  \DeclareOldFontCommand{\bf}{\normalfont\bfseries}{\mathbf}
  \DeclareOldFontCommand{\it}{\normalfont\itshape}{\mathit}
  \DeclareOldFontCommand{\sl}{\normalfont\slshape}{\@nomath\sl}
  \DeclareOldFontCommand{\sc}{\normalfont\scshape}{\@nomath\sc}
\fi

\def\alpha{{\Greekmath 010B}}%
\def\beta{{\Greekmath 010C}}%
\def\gamma{{\Greekmath 010D}}%
\def\delta{{\Greekmath 010E}}%
\def\epsilon{{\Greekmath 010F}}%
\def\zeta{{\Greekmath 0110}}%
\def\eta{{\Greekmath 0111}}%
\def\theta{{\Greekmath 0112}}%
\def\iota{{\Greekmath 0113}}%
\def\kappa{{\Greekmath 0114}}%
\def\lambda{{\Greekmath 0115}}%
\def\mu{{\Greekmath 0116}}%
\def\nu{{\Greekmath 0117}}%
\def\xi{{\Greekmath 0118}}%
\def\pi{{\Greekmath 0119}}%
\def\rho{{\Greekmath 011A}}%
\def\sigma{{\Greekmath 011B}}%
\def\tau{{\Greekmath 011C}}%
\def\upsilon{{\Greekmath 011D}}%
\def\phi{{\Greekmath 011E}}%
\def\chi{{\Greekmath 011F}}%
\def\psi{{\Greekmath 0120}}%
\def\omega{{\Greekmath 0121}}%
\def\varepsilon{{\Greekmath 0122}}%
\def\vartheta{{\Greekmath 0123}}%
\def\varpi{{\Greekmath 0124}}%
\def\varrho{{\Greekmath 0125}}%
\def\varsigma{{\Greekmath 0126}}%
\def\varphi{{\Greekmath 0127}}%

\def\nabla{{\Greekmath 0272}}
\def\FindBoldGroup{%
   {\setbox0=\hbox{$\mathbf{x\global\edef\theboldgroup{\the\mathgroup}}$}}%
}

\def\Greekmath#1#2#3#4{%
    \if@compatibility
        \ifnum\mathgroup=\symbold
           \mathchoice{\mbox{\boldmath$\displaystyle\mathchar"#1#2#3#4$}}%
                      {\mbox{\boldmath$\textstyle\mathchar"#1#2#3#4$}}%
                      {\mbox{\boldmath$\scriptstyle\mathchar"#1#2#3#4$}}%
                      {\mbox{\boldmath$\scriptscriptstyle\mathchar"#1#2#3#4$}}%
        \else
           \mathchar"#1#2#3#4%
        \fi
    \else
        \FindBoldGroup
        \ifnum\mathgroup=\theboldgroup 
           \mathchoice{\mbox{\boldmath$\displaystyle\mathchar"#1#2#3#4$}}%
                      {\mbox{\boldmath$\textstyle\mathchar"#1#2#3#4$}}%
                      {\mbox{\boldmath$\scriptstyle\mathchar"#1#2#3#4$}}%
                      {\mbox{\boldmath$\scriptscriptstyle\mathchar"#1#2#3#4$}}%
        \else
           \mathchar"#1#2#3#4%
        \fi     	
	  \fi}

\newif\ifGreekBold  \GreekBoldfalse
\let\SAVEPBF=\pbf
\def\pbf{\GreekBoldtrue\SAVEPBF}%

\@ifundefined{theorem}{\newtheorem{theorem}{Theorem}}{}
\@ifundefined{lemma}{\newtheorem{lemma}[theorem]{Lemma}}{}
\@ifundefined{corollary}{\newtheorem{corollary}[theorem]{Corollary}}{}
\@ifundefined{conjecture}{}{}
\@ifundefined{proposition}{\newtheorem{proposition}[theorem]{Proposition}}{}
\@ifundefined{axiom}{}{}
\@ifundefined{remark}{\newtheorem{remark}{Remark}}{}
\@ifundefined{example}{}{}
\@ifundefined{exercise}{}{}
\@ifundefined{definition}{\newtheorem{definition}{Definition}}{}

\@ifundefined{mathletters}{%
  \newcounter{equationnumber}
  \def\mathletters{%
     \addtocounter{equation}{1}
     \edef\@currentlabel{\theequation}%
     \setcounter{equationnumber}{\c@equation}
     \setcounter{equation}{0}%
     \edef\theequation{\@currentlabel\noexpand\alph{equation}}%
  }
  
}{}

\@ifundefined{BibTeX}{%
    \def\BibTeX{{\rm B\kern-.05em{\sc i\kern-.025em b}\kern-.08em
                 T\kern-.1667em\lower.7ex\hbox{E}\kern-.125emX}}}{}%
\@ifundefined{AmS}%
    {\def\AmS{{\protect\usefont{OMS}{cmsy}{m}{n}%
                A\kern-.1667em\lower.5ex\hbox{M}\kern-.125emS}}}{}%
\@ifundefined{AmSTeX}{}{}%
\def\@@eqncr{\let\@tempa\relax
    \ifcase\@eqcnt \def\@tempa{& & &}\or \def\@tempa{& &}%
      \else \def\@tempa{&}\fi
     \@tempa
     \if@eqnsw
        \iftag@
           \@taggnum
        \else
           \@eqnnum\stepcounter{equation}%
        \fi
     \fi
     \global\tag@false
     \global\@eqnswtrue
     \global\@eqcnt\z@\cr}

\def\TCItag{\@ifnextchar*{\@TCItagstar}{\@TCItag}}
\def\@TCItag#1{%
    \global\tag@true
    \global\def\@taggnum{(#1)}%
    \global\def\@currentlabel{#1}}
\def\@TCItagstar*#1{%
    \global\tag@true
    \global\def\@taggnum{#1}%
    \global\def\@currentlabel{#1}}
\def\tint{\msi@int\textstyle\int}%
\def\tiint{\msi@int\textstyle\iint}%
\def\tiiint{\msi@int\textstyle\iiint}%
\def\tiiiint{\msi@int\textstyle\iiiint}%
\def\tidotsint{\msi@int\textstyle\idotsint}%
\def\toint{\msi@int\textstyle\oint}%

\def\tprod{\mathop{\textstyle \prod }}%
\newtoks\temptoksa
\newtoks\temptoksb
\newtoks\temptoksc

\def\msi@int#1#2{%
 \def\@temp{{#1#2\the\temptoksc_{\the\temptoksa}^{\the\temptoksb}}}%
 \futurelet\@nextcs
 \@int
}

\def\@int{%
   \ifx\@nextcs\limits
      \typeout{Found limits}%
      \temptoksc={\limits}%
	  \let\@next\@intgobble%
   \else\ifx\@nextcs\nolimits
      \typeout{Found nolimits}%
      \temptoksc={\nolimits}%
	  \let\@next\@intgobble%
   \else
      \typeout{Did not find limits or no limits}%
      \temptoksc={}%
      \let\@next\msi@limits%
   \fi\fi
   \@next
}%

\def\@intgobble#1{%
   \typeout{arg is #1}%
   \msi@limits
}

\def\msi@limits{%
   \temptoksa={}%
   \temptoksb={}%
   \@ifnextchar_{\@limitsa}{\@limitsb}%
}

\def\@limitsa_#1{%
   \temptoksa={#1}%
   \@ifnextchar^{\@limitsc}{\@temp}%
}

\def\@limitsb{%
   \@ifnextchar^{\@limitsc}{\@temp}%
}

\def\@limitsc^#1{%
   \temptoksb={#1}%
   \@ifnextchar_{\@limitsd}{\@temp}%
}

\def\@limitsd_#1{%
   \temptoksa={#1}%
   \@temp
}

\def\dint{\msi@int\displaystyle\int}%
\def\diint{\msi@int\displaystyle\iint}%
\def\diiint{\msi@int\displaystyle\iiint}%
\def\diiiint{\msi@int\displaystyle\iiiint}%
\def\didotsint{\msi@int\displaystyle\idotsint}%
\def\doint{\msi@int\displaystyle\oint}%

\begin{document}

\title{Primal and Dual Characterizations for Farkas type Lemmas in Terms of Closedness Criteria}

\author*[1]{\fnm{Dinh} \sur{N.}}\email{ndinh@hcmiu.edu.vn}
\author[2]{\fnm{Goberna}\sur{M.A.}}\email{mgoberna@ua.es}
\author[3]{\fnm{Volle}\sur{M.}}\email{michel.volle@univ-avignon.fr}

\affil*[1]{\orgdiv{Department of Mathematics, International University, Vietnam National, University-HCMC}, \city{Ho Chi Minh city}, \country{Vietnam}}
\affil[2]{\orgdiv{Department of Mathematics, University of Alicante}, \postcode{03080}, \city{Alicante}, \country{Spain}}
\affil[3]{\orgdiv{D\'{e}partement de Math\'{e}matiques, LMA EA 2151},  \city{Avignon}, \country{France}}

\maketitle
\date{}

\begin{abstract}
_This paper deals with the characterization, in terms of closedness of
certain sets regarding other sets, of Farkas lemmas determining when the
upperlevel set of a convex function $f$ contains a set of the form $C\cap
\mathbb{A}^{-1}\left( D\right) ,$ where $C$ and $D$ are convex sets (not
necessarily cones) in locally convex spaces $X$ (with topological dual $%
X^{\prime }$) and $Y,$ respectively, while $\mathbb{A}$ is a continuous
linear operator from $X$ to $Y$. More in detail, each of the mentioned
characterizations of Farkas type lemmas consists in the closedness of
certain subset of either one of the "primal" spaces $X\times Y\times \mathbb{%
R}$ and $Y\times \mathbb{R},$ or of the "dual" space $X^{\prime }\times
\mathbb{R}$, regarding some singleton set of the corresponding space.
Moreover, the paper also provides an existence theorem for the feasible set $%
C\cap \mathbb{A}^{-1}\left( D\right) $ in terms of the closedness of certain
subset of the dual space $X^{\prime }\times \mathbb{R}$ regarding the
singleton set formed by the null element. These results are illustrated with
significant applications.
\end{abstract}



\section{Introduction}

Let us start with an informal preamble. In the framework of locally convex
Hausdorff topological vector spaces (lcHtvs in brief), the homogeneous
Farkas' lemma says that the statements

\textrm{(i)} $\mathbb{A}x\in Q\Longrightarrow \left\langle \QTR{frametitle}{a%
}^{\prime },x\right\rangle \geq 0$

and

\textrm{(ii)} $\exists \mu \in Q^{+}$ such that $\mathbb{A}^{\#}\mu =%
\QTR{frametitle}{a}^{\prime }$\newline
are equivalent if and only if (see Corollary \ref{Cor6} below) $\mathbb{A}%
^{\#}\left( Q^{+}\right) $ is $w^{\ast }$ $-$ closed regarding $\left\{
\QTR{frametitle}{a}^{\prime }\right\} ,$ that means
\begin{equation*}
\left\{ \QTR{frametitle}{a}^{\prime }\right\} \cap w^{\ast }-\limfunc{cl}%
\mathbb{A}^{\#}\left( Q^{+}\right) =\left\{ \QTR{frametitle}{a}^{\prime
}\right\} \cap \mathbb{A}^{\#}\left( Q^{+}\right)
\end{equation*}%
or, in other words,%
\begin{equation*}
\text{either }\QTR{frametitle}{a}^{\prime }\in \mathbb{A}^{\#}\left(
Q^{+}\right) \text{ or }\QTR{frametitle}{a}^{\prime }\notin w^{\ast }-%
\limfunc{cl}\mathbb{A}^{\#}\left( Q^{+}\right) .
\end{equation*}%
Here, $\mathbb{A}:X\longrightarrow Y$ is a continuous linear operator, $%
\mathbb{A}^{\#}$ is its adjoint, $Q$ is a closed convex cone in $Y,$ whose
topological dual is $Y^{\prime }$ , $Q^{+}:=\left\{ y^{\prime }\in Y^{\prime
}:\left\langle \QTR{frametitle}{y^{\prime }},y\right\rangle \geq 0,\forall
\QTR{frametitle}{y}\in Q\right\} $ is the positive dual cone of $Q,$ and $%
\QTR{frametitle}{a}^{\prime }$ is a continuous linear form on $X.$

It is less known that, besides this $w^{\ast }-$closedness criterion, there
is another one formulated in the primal product space $Y\times \mathbb{R}$
giving another characterization, namely that \textrm{(i)}$%
\Longleftrightarrow $\textrm{(ii)} if and only if (see Corollary\ \ref{Cor2}
below)%
\begin{equation*}
\bigcup\limits_{x\in X}\left\{ \left( \mathbb{A}x,\left\langle
\QTR{frametitle}{a}^{\prime },x\right\rangle \right) \right\} +\left(
-Q\right) \times \mathbb{R}_{+}
\end{equation*}%
is closed regarding $\left\{ \left( 0_{Y},-1\right) \right\} ,$ where $0_{Y}$
denotes the null element of $Y.$

Notice that, for this second criterion, the convex cone $Q$ is not
necessarily closed and does not necessarily contains the origin. The aim of
this note is to formulate similar results for general Farkas type lemmas
with a special attention concerning the infinite nonhomogeneous linear
systems.

More formally, given two lcHtvs $X$ and $Y$ with respective topological
duals $X^{\prime }$ and $Y^{\prime },$ two nonempty convex sets $C$ and $D$
of $X$ and $Y,$ respectively, and an extended real-valued proper convex
function $f:X\longrightarrow \overline{\mathbb{R}}:=\mathbb{R\cup }\left\{
\pm \infty \right\} ,$ we are concerned with the property

$(\mathcal{A})$\ \ \ \ $x \in C \cap \mathbb{A}^{-1} (D) \ \ \Longrightarrow
f(x) \geq 0$,

\noindent which is in particular satisfied (see Theorem \ref{Th3}) under the
dual condition

$(\mathcal{B}) $ \ \ $\left\{
\begin{array}{ll}
& \exists \left( u^{\prime },v^{\prime },\lambda \right) \in \func{dom}%
f^{\ast }\times \beta \left( C\right) \times \beta \left( D\right) \ \text{
such\ that } \\
& f^{\ast }\left( u^{\prime }\right) +\sigma _{C}\left( v^{\prime }\right)
+\sigma _{D}\left( \lambda \right) \leq 0\text{ and }u^{\prime }+v^{\prime
}=-\mathbb{A}^{\#}\lambda .%
\end{array}
\right. $

\noindent Here $f^{\ast }$ denotes the conjugate of $f,$ i.e., $f^{\ast
}\left( x^{\prime }\right) =\sup_{x\in X}\left( \left\langle x^{\prime
},x\right\rangle -f\left( x\right) \right) ,$ $\sigma _{C}$ the support
function of $C,$ i.e. the conjugate of the indicator $\delta _{C}$ of $C$
(with $\delta _{C}\left( x\right) =0$ if $x\in C$ and $\delta _{C}\left(
x\right) =+\infty $ otherwise), and $\beta \left( C\right) $ the barrier
cone of $C,$ i.e., the domain $\func{dom}\sigma _{C}=\left\{ x^{\prime }\in
X^{\prime }:\sigma _{C}\left( x^{\prime }\right) <+\infty \right\} $ of $%
\sigma _{C}.$ Our purpose is to characterize the direct implication, that
is, $(\mathcal{A}) \Longrightarrow (\mathcal{B})$ or, equivalently, $(%
\mathcal{A}) \Longleftrightarrow (\mathcal{B})$, in terms of a primal (resp.
dual) closedness criterion. To this aim, let us recall the concept of
closedness regarding a set \cite{Bot10}.

Given two subsets $U$ and $V$ of a topological space $Z,$ one says that $U$
is closed regarding $V$ if $\overline{U}\cap V=U\cap V,$ where $\overline{U}$
is the closure of $U$ in $Z.$ In particular, Given $z\in Z,$ $U$ is closed
regarding $\left\{ z\right\} $ if and only if either $z\in U$ or $z\notin
\overline{U},$ $U$ is closed regarding $V$ if and only if $U$ is closed
regarding $\left\{ z\right\} $ for all $z\in V,$ and $U$ is closed if and
only if it is closed regarding $Z.$

The term \textit{characterization of Farkas' lemma} was introduced by V.
Jeyakumar, S. Kum, and G. M. Lee \cite{JKL08} and that of \textit{%
characterization of stable Farkas' lemma} by V. Jeyakumar and G. M. Lee \cite%
{JL08}, both papers published in 2008. Their characterization of stable
Farkas' lemma consisted in the closedness of certain cone associated with
the given conic system. After the introduction of the concept of closedness
regarding a given set by Bo\c{t} in 2010, successive characterizations of
(stable) Farkas' lemma in terms of closedness of certain set regarding
another set were obtained for different types of systems, e.g., systems
involving convex, nonconvex composite functions \cite{DMVV17}, \cite{DVV14},
\cite{DVN-08}, systems involving vector-valued functions (\cite{DGLL19} and
\cite{DGLM17}), systems involving a family of finite subsets of the given
index set \cite{DGLV22}, etc. The characterization of stable Farkas' lemma
for\ nonconvex composite semi-infinite programming has been considered in
\cite{LZ15} and \cite{LPW19}.

Considering feasible sets of the form $C\cap \mathbb{A}^{-1}\left( D\right)
, $ with the convex set $D$ not being necessarily a cone, is the main
novelty of this paper, which is organized as follows. Section 2
characterizes the equivalence $(\mathcal{A})\Longleftrightarrow (\mathcal{B}%
) $ in terms of the closedness regarding $\left\{ \left(
0_{X},0_{Y},-1\right) \right\} $ and $\left\{ \left( 0_{Y},-1\right)
\right\} $\ of certain subsets of the primal spaces $X\times Y\times \mathbb{%
R}$ and $Y\times \mathbb{R},$ respectively (see Theorems \ref{Th1} and \ref%
{Th2}). Section 3, in turn, characterizes the equivalence $(\mathcal{A}%
)\Longleftrightarrow (\mathcal{B})$ in terms of the closedness regarding $%
\left\{ \left( 0_{X^{\prime }},0\right) \right\} $ of certain subset of the
dual space $X^{\prime }\times \mathbb{R}$ (Theorem \ref{Th3}) as well as the
non-emptyness of the feasible set $C\cap \mathbb{A}^{-1}\left( D\right) $
(Proposition \ref{prop1}) and the stable Farkas' lemma (Proposition \ref%
{prop2}). Section 4 is devoted to linear infinite systems. Section 5
characterizes Farkas-type lemmas oriented to identify the minima of convex
(resp. concave) functions on $C\cap \mathbb{A}^{-1}\left( D\right) ,$ reason
why they are known in the literature as convex (resp. concave) Farkas'
lemmas. More in detail, the convex Farkas' lemmas provided in Section 3
determine when a convex feasible set of the form $C\cap \mathbb{A}%
^{-1}\left( D\right) $ is contained in the reverse-convex set $\left\{ x\in
X:f\left( x\right) \geq 0\right\} ,$ where $f$ is convex. Similarly, Section
5 characterizes the containment of $C\cap \mathbb{A}^{-1}\left( D\right) $
in the sublevel set of a convex function $f$ \ (Theorem \ref{thmV1}).
Finally, Section 6 shows applications to constrained convex minimization
problems (optimization, strong duality and stable strong duality theorems)
and to functional approximation\ by polynomials (Farkas-type lemmas).\bigskip

\section{Primal closedness characterizations}

Recall that $f:X\longrightarrow \overline{\mathbb{R}}$ ($f\in \overline{%
\mathbb{R}}^{X}$ in short) is a proper convex function. We consider the
vector-valued mapping $H:X\times Y\times \func{dom}f\longrightarrow X\times
Y\times \mathbb{R}$ such that
\begin{equation*}
H\left( x,y,v\right) =\left( x-v,\mathbb{A}x-y,f\left( v\right) \right) ,
\end{equation*}%
which is convex with respect to the cone $\left\{ 0_{X}\right\} \times
\left\{ 0_{Y}\right\} \times \mathbb{R}_{+}.$ Therefore,%
\begin{equation*}
\begin{array}{ll}
\mathcal{F} & :=H\left( C\times D\times \func{dom}f\right) +\left\{
0_{X}\right\} \times \left\{ 0_{Y}\right\} \times \mathbb{R}_{+}\medskip \\
& =\bigcup\limits_{x\in C}\text{ }\left\{ \left( x,\mathbb{A}x,0\right)
\right\} +\bigcup\limits_{v\in \func{dom}f}\left\{ \left( -v,0_{Y},f\left(
v\right) \right) \right\} +\left\{ 0_{X}\right\} \times \left( -D\right)
\times \mathbb{R}_{+}%
\end{array}%
\end{equation*}%
is a nonempty convex subset of $X\times Y\times \mathbb{R}$ and $\mathbb{R}%
_{+}\mathcal{F}$ is a convex cone in $X\times Y\times \mathbb{R}.\medskip $

\begin{theorem}[1st characterization of Farkas' Lemma]
\label{Th1}Let $f\in \overline{\mathbb{R}}^{X}$ be proper convex, $C$ (resp.
$D$) be a nonempty convex subset of $X$ (resp. $Y$). Then, one has $(%
\mathcal{A})\Longrightarrow (\mathcal{B})$ (hence, $(\mathcal{A}%
)\Longleftrightarrow (\mathcal{B})$) if and only if $\mathbb{R}_{+}\mathcal{F%
}$ is closed regarding $\left\{ \left( 0_{X},0_{Y},-1\right) \right\} .$
\end{theorem}

\medskip \noindent \textbf{Proof.} \textit{Sufficiency}\textbf{.} Assume
that $(\mathcal{A})$ holds. We first prove that $\left(
0_{X},0_{Y},-1\right) \notin \mathbb{R}_{+}\mathcal{F}$. Otherwise, by the
definition of $\mathcal{F},$ there exist $\theta \in \mathbb{R}_{+}$ and $%
\left( x,v,d,\eta \right) \in C\times \func{dom}f\times D\times \mathbb{R}%
_{+}$ such that
\begin{equation*}
\left( 0_{X},0_{Y},-1\right) =\theta \left[ \left( x,\mathbb{A}x,0\right)
+\left( -v,0_{Y},f\left( v\right) \right) +\left( 0_{X},-d,\eta \right) %
\right] .
\end{equation*}%
Clearly, $\theta >0,$ $0_{X}=x-v,$ $0_{Y}=\mathbb{A}x-d,$ and $f\left(
v\right) +\eta <0.$ Therefore, $x=v\in C,$ $\mathbb{A}x=d\in D$ and $f\left(
v\right) <-\eta \leq 0.$ Consequently, $x\in C\cap \mathbb{A}^{-1}\left(
D\right) $ with $f\left( x\right) <0$ in contradiction with $(\mathcal{A})$.

Since $\mathbb{R}_{+}\mathcal{F}$ is closed regarding $\left\{ \left(
0_{X},0_{Y},-1\right) \right\} $ and $\left( 0_{X},0_{Y},-1\right) \notin
\mathbb{R}_{+}\mathcal{F}$, we have $\left( 0_{X},0_{Y},-1\right) \notin
\overline{\mathbb{R}_{+}\mathcal{F}},$ which is a nonempty closed convex
cone. Applying the Hahn-Banach Theorem, there exists a non-zero vector $%
\left( x^{\prime },\mu ,s\right) \in X^{\prime }\times Y^{\prime }\times
\mathbb{R}$ such that%
\begin{equation*}
\left\langle x^{\prime },0_{X}\right\rangle +\left\langle \mu
,0_{Y}\right\rangle -s=-s>\gamma ,
\end{equation*}%
where
\begin{equation*}
\gamma :=\sup\limits_{\left( x,v,d,\left( \theta ,\alpha \right) \right) \in
C\times \func{dom}f\times D\times \mathbb{R}^{2}}\theta \left( \left\langle
x^{\prime },x-v\right\rangle +\left\langle \mu ,\mathbb{A}x-d\right\rangle
+s\left( f\left( v\right) +\alpha \right) \right) .
\end{equation*}
We have necessarily $\gamma =0,$ $s<0$ and%
\begin{equation*}
\left\langle x^{\prime },x-v\right\rangle +\left\langle \mu ,\mathbb{A}%
x-d\right\rangle +sf\left( v\right) \leq 0,\forall \left( x,v,d\right) \in
C\times \func{dom}f\times D.
\end{equation*}%
So,%
\begin{equation*}
\sup\limits_{x\in C}\left\langle x^{\prime }+\mathbb{A}^{\#}\mu
,x\right\rangle +\sup\limits_{v\in \func{dom}f}\left( -\left\langle
x^{\prime },v\right\rangle +sf\left( v\right) \right) +\sup\limits_{d\in
D}\left\langle \mu ,-d\right\rangle \leq 0,
\end{equation*}%
that is,%
\begin{equation*}
\sigma _{C}\left( x^{\prime }+\mathbb{A}^{\#}\mu \right) -sf^{\ast }\left(
\frac{x^{\prime }}{s}\right) +\sigma _{D}\left( -\mu \right) \leq 0,
\end{equation*}%
or, equivalently,%
\begin{equation*}
f^{\ast }\left( \frac{x^{\prime }}{s}\right) +\sigma _{C}\left( -\frac{%
x^{\prime }+\mathbb{A}^{\#}\mu }{s}\right) +\sigma _{D}\left( \frac{\mu }{s}%
\right) \leq 0.
\end{equation*}
Setting $\left( \mu ^{\prime },\lambda \right) :=s^{-1}\left( x^{\prime
},\mu \right) $ we obtain that $(\mathcal{B})$ holds.

{\textit{Necessity}}. Assume $(\mathcal{A}) \Longrightarrow (\mathcal{B})$
and suppose that $\left( 0_{X},0_{Y},-1\right) \notin \mathbb{R}_{+}\mathcal{%
F}.$ We must prove that $\left( 0_{X},0_{Y},-1\right) \notin \overline{%
\mathbb{R}_{+}\mathcal{F}}.$

Since $\left( 0_{X},0_{Y},-1\right) \notin \mathbb{R}_{+}\mathcal{F},$
statement $(\mathcal{A})$ holds. Otherwise there exist $\overline{x}\in C,$ $%
\overline{d}\in D,$ $\mathbb{A}\overline{x}=\overline{d},$ and $f\left(
\overline{x}\right) <0.$ By definition of $\mathcal{F}$ we then have%
\begin{equation*}
\left( \overline{x},\mathbb{A}\overline{x},0\right) +\left( -\overline{x}%
,0_{Y},f\left( \overline{x}\right) \right) +\left( 0_{X},-\overline{d}%
,0\right) =\left( 0_{X},0_{Y},f\left( \overline{x}\right) \right) \in
\mathcal{F}.
\end{equation*}%
We have $\left( 0_{X},0_{Y},-1\right) \notin \left( \mathbb{R}_{+}\diagdown
\left\{ 0\right\} \right) \mathcal{F}$ in contradiction with $\left(
0_{X},0_{Y},-1\right) \notin \mathbb{R}_{+}\mathcal{F}.$

Therefore, statement $(\mathcal{A})$ holds and, by hypothesis, $( \mathcal{B}%
)$ holds too.

For all $\left( x,v,d,\eta \right) \in C\times \func{dom}f\times D\times
\mathbb{R}_{+}$ we then have%
\begin{equation*}
\left\langle u^{\prime },v\right\rangle -f\left( v\right) +\left\langle
v^{\prime },x\right\rangle +\left\langle \lambda ,d\right\rangle \leq
f^{\ast }\left( u^{\prime }\right) +\sigma _{C}\left( v^{\prime }\right)
+\sigma _{D}\left( \lambda \right) \leq 0.
\end{equation*}

Since $u^{\prime }+v^{\prime }=-\mathbb{A}^{\#}\lambda ,$ it follows that%
\begin{equation}
-\left\langle u^{\prime },x-v\right\rangle -\left\langle \lambda ,\mathbb{A}%
x-d\right\rangle \leq f\left( v\right) +\eta .  \label{1}
\end{equation}

Consider the continuous linear form
\begin{equation*}
\xi \left( u,y\right) :=-\left\langle u^{\prime },u\right\rangle
-\left\langle \lambda ,y\right\rangle ,\forall \left( u,y\right) \in X\times
Y,
\end{equation*}%
and let $\left( u,y,r\right) \in \mathcal{F}.$ There exists $\left(
x,v,d,\eta \right) \in C\times \func{dom}f\times D\times \mathbb{R}_{+}$
such that
\begin{equation*}
\left( u,y,r\right) =\left( x-v,\mathbb{A}x-d,f\left( v\right) +\eta \right)
.
\end{equation*}%
By (\ref{1}), $\left( u,y,r\right) $ satisfies $\xi \left( u,y\right) \leq
r. $ So, $\mathcal{F}$ is contained in the closed convex cone epigraph of $%
\xi , $ denoted by $\limfunc{epi}\xi .$ Consequently, the closed convex cone
$\overline{\mathbb{R}_{+}\mathcal{F}}$ is contained in $\limfunc{epi}\xi $.
Since $\left( 0_{X},0_{Y},-1\right) \notin \limfunc{epi}\xi ,$ we conclude
that $\left( 0_{X},0_{Y},-1\right) \notin \overline{\mathbb{R}_{+}\mathcal{F}%
}$ and we are done. \hfill $\square \medskip \medskip $

We now associate a perturbation function with statements $(\mathcal{A})$ and
$(\mathcal{B})$, namely $F:X\times \left( X\times Y\right) \longrightarrow
\overline{\mathbb{R}}$ such that%
\begin{equation}  \label{eq22nw}
F\left( x,\left( u,y\right) \right) :=\left\{
\begin{array}{ll}
f\left( x-u\right) , & \text{if }x\in C, u \in X, \text{ and }\mathbb{A}x\in
y+D, \\
+\infty , & \text{else.}%
\end{array}%
\right.
\end{equation}

We may rewrite $(\mathcal{A})$ and $(\mathcal{B})$ in terms of $F$\ as
follows:

$(\mathcal{A})$\ \ \ \ \ $F\left( x,\left( 0_{X},0_{Y}\right) \right) \geq
0,\forall x\in X,$

$(\mathcal{B}) $ \ \ $\left\{
\begin{array}{ll}
& \exists \left( u^{\prime },\lambda \right) \in X^{\prime }\times Y^{\prime
} \text{ such \ that} \\
& F\left( x,\left( u,y\right) \right) +\left\langle u^{\prime
},u\right\rangle +\left\langle \lambda ,y\right\rangle \geq 0,\forall \left(
x,\left( u,y\right) \right) \in X\times \left( X\times Y\right) .%
\end{array}
\right. $

The $S_{F}-$procedure introduced in \cite[Theorem 2.2]{VB2024} in a more
general nonconvex setting is said to be valid if $(\mathcal{A}%
)\Longrightarrow (\mathcal{B})$ or, equivalently, if $(\mathcal{A}%
)\Longleftrightarrow (\mathcal{B})$. Note that the set $\mathcal{F}$
coincides with the projection on $\left( X\times Y\right) \times \mathbb{R}$
of the epigraph of $F,$ which is convex. Consequently, \cite[Theorem 2.2]%
{VB2024} can be applied to the present situation replacing the closed convex
hull of $\mathbb{R}_{+}\mathcal{F}$ by its closure.

Let us now consider the more simple perturbation function%
\begin{equation*}
F_{0}\left( x,y\right) :=\left\{
\begin{array}{ll}
f\left( x\right) , & \text{if }x\in C\text{ and }\mathbb{A}x\in y+D, \\
+\infty , & \text{else,}%
\end{array}%
\right.
\end{equation*}%
which is convex (and proper if $C\cap \func{dom}f\neq \emptyset $). Here
again $(\mathcal{A})$ reads%
\begin{equation*}
F_{0}\left( x,0_{Y}\right) \geq 0,\forall x\in X,
\end{equation*}%
and the $S_{F_{0}}-$procedure is valid if and only if statement $(\mathcal{A}%
)$ entails that%
\begin{equation*}
\exists \lambda \in \beta \left( D\right) \text{ such that }F_{0}\left(
x,y\right) +\left\langle \lambda ,y\right\rangle \geq 0,\forall \left(
x,y\right) \in X\times Y,
\end{equation*}%
that is to say
\begin{equation}
\exists \lambda \in \beta \left( D\right) \text{ such that }f\left( x\right)
+\left\langle \lambda ,\mathbb{A}x\right\rangle \geq \sigma _{D}\left(
\lambda \right) ,\forall x\in C.  \label{2.0}
\end{equation}

The projection of $\limfunc{epi}F_{0}$ on $Y\times \mathbb{R}$ is the
nonempty convex set%
\begin{equation*}
\mathcal{F}_{0}:=\bigcup\limits_{x\in C\cap \func{dom}f}\text{ }\left\{
\left( \mathbb{A}x,f\left( x\right) \right) \right\} +\left( -D\right)
\times \mathbb{R}_{+}.
\end{equation*}%
Applying \cite[Theorem 2.2]{VB2024} we obtain\medskip

\begin{theorem}[2nd characterization of Farkas' Lemma]
\label{Th2}Let $C$ (resp. $D$) be a nonempty convex subset of $X$ (resp. $Y$%
) and $f\in \overline{\mathbb{R}}^{X}$ be convex proper and such that $C\cap
\func{dom}f\neq \emptyset .$ Then, the statement

$(\mathcal{A})$\ \ \ \ $x\in C\cap \mathbb{A}^{-1}\left( D\right)
\Longrightarrow f\left( x\right) \geq 0$

\noindent is equivalent to

$(\mathcal{B}_0) $ \ \ (\ref{2.0}) holds,\newline
if and only if the convex cone%
\begin{equation*}
\mathbb{R}_{+}\mathcal{F}_{0}=\mathbb{R}_{+}\left[ \bigcup\limits_{x\in
C\cap \func{dom}f}\text{ }\left\{ \left( \mathbb{A}x,f\left( x\right)
\right) \right\} +\left( -D\right) \times \mathbb{R}_{+}\right] \text{ }%
\medskip
\end{equation*}%
is closed regarding $\left\{ \left( 0_{Y},-1\right) \right\}.$
\end{theorem}

\begin{corollary}
\label{Cor1}Let $C\subset X$ and $D\subset Y$ be nonempty convex sets. Given
$\left( x^{\prime },r\right) \in X^{\prime }\times \mathbb{R},$ the
statements

$(\mathcal{A}_{1})$\ \ \ \ $x\in C\cap \mathbb{A}^{-1}\left( D\right)
\Longrightarrow \left\langle x^{\prime },x\right\rangle \leq r$\newline
and

$(\mathcal{B}_{1})$ $\exists \lambda \in \beta \left( D\right) $ such that $%
\sigma _{C}\left( x^{\prime }-\mathbb{A}^{\#}\lambda \right) +\sigma
_{D}\left( \lambda \right) \leq r$

\noindent are equivalent if and only if the set\textrm{\ }%
\begin{equation*}
\mathbb{R}_{+}\left[ \bigcup\limits_{x\in C\cap \func{dom}f}\text{ }\left\{
\left( \mathbb{A}x,-\left\langle x^{\prime },x\right\rangle \right) \right\}
+\left\{ \left( 0_{Y},r\right) \right\} +\left( -D\right) \times \mathbb{R}%
_{+}\right] \text{ }
\end{equation*}%
is closed regarding $\left\{ \left( 0_{Y},-1\right) \right\} .$
\end{corollary}

\textbf{Proof.} Apply Theorem \ref{Th2} to $f\left( x\right) =r-\left\langle
x^{\prime },x\right\rangle .$\hfill $\square \medskip $

\begin{corollary}
\label{Cor2}Let $a^{\prime }\in X^{\prime }$ and $Q$ be a nonempty convex
cone in $Y.$ Then, the statements $(\mathcal{A}_{2})$ and $(\mathcal{B}_{2})$
below,

$(\mathcal{A}_{2})$ \ \ \ $\mathbb{A}x\in Q\Longrightarrow \left\langle
a^{\prime },x\right\rangle \geq 0$\newline
and

$(\mathcal{B}_{2})$ \ \ $\exists \mu \in Q^{+}$ such that $\mathbb{A}%
^{\#}\mu =a^{\prime },$ \newline
are equivalent if and only if the set%
\begin{equation*}
\bigcup\limits_{x\in X}\text{ }\left\{ \left( \mathbb{A}x,\left\langle
a^{\prime },x\right\rangle \right) \right\} +\left( -Q\right) \times \mathbb{%
R}_{+}\text{ is closed regarding }\left\{ \left( 0_{Y},-1\right) \right\} .
\end{equation*}
\end{corollary}

\textbf{Proof.} Apply Corollary \ref{Cor1} with $D=Q,$ $\left( x^{\prime
},r\right) =\left( -a^{\prime },0\right) ,$ $C=X,$ $\sigma _{C}=\delta
_{\left\{ 0_{X^{\prime }}\right\} },$ and $\sigma _{D}=\delta _{Q^{-}},$
where $Q^{-}:=-Q^{+}$ is the negative dual cone of $Q.$\ Then, statement $(%
\mathcal{B}_1)$ reads as: there exists $\lambda \in Q^{-}$ such that $%
\mathbb{A}^{\#}\lambda =-a^{\prime },$ that is, $(\mathcal{B}_2)$ holds with
$\mu =-\lambda \in Q^{-}.$ To conclude the proof, let us check the two
closedness conditions in Corollaries \ref{Cor1} and \ref{Cor2} are
equivalent for $\left( x^{\prime },r\right) =\left( -a^{\prime },0\right) ,$
which is clear if the convex cone $Q$ contains $0_{Y},$ but this condition
is not required. One has in fact%
\begin{equation*}
\begin{array}{c}
\mathbb{R}_{+}\!\!\!\left[\!\! \bigcup\limits_{x\in X}\!\!\left\{ \left(
\mathbb{A}x,\left\langle a^{\prime },x\right\rangle \right) \right\} +\left(
-Q\right) \times \mathbb{R}_{+}\right] \!\!=\!\!\left\{ \left(
0_{Y},0\right) \right\}\! \cup\! \left[\! \bigcup\limits_{x\in X}\!\!\left\{
\left( \mathbb{A}x,\left\langle a^{\prime },x\right\rangle \right) \right\}
+\left( -Q\right) \times \mathbb{R}_{+}\!\!\right] .%
\end{array}%
\end{equation*}%
Thus,%
\begin{equation*}
\overline{\mathbb{R}_{+}\left[ \bigcup\limits_{x\in X}\left\{ \left( \mathbb{%
A}x,\left\langle a^{\prime },x\right\rangle \right) \right\} +\left(
-Q\right) \times \mathbb{R}_{+}\right] }=\overline{\bigcup\limits_{x\in
X}\left\{ \left( \mathbb{A}x,\left\langle a^{\prime },x\right\rangle \right)
\right\} +\left( -Q\right) \times \mathbb{R}_{+}}
\end{equation*}%
and, by this, the two closedness statements in Corollaries \ref{Cor1} and %
\ref{Cor2} coincide when $\left( x^{\prime },r\right) =\left( -a^{\prime
},0\right) .$\hfill $\square \medskip $

\section{Dual closedness characterization}

We consider the topological dual $X^{\prime }$ of $X$ equipped with the $%
w^{\ast }$ topology and $X^{\prime }\times \mathbb{R}$ provided with the
product topology. The closure of a set $A\subset X^{\prime }\times \mathbb{R}
$ is denoted by either $w^{\ast }-\limfunc{cl}A$ or by $\overline{A}$ (in
short).

We associate with $D\subset Y,$ $D\neq \emptyset ,$ the sets%
\begin{equation}
B:=\mathbb{A}^{-1}\left( D\right) \text{ and }\Lambda
:=\bigcup\limits_{\lambda \in \beta \left( D\right) }\left\{ \left( \mathbb{A%
}^{\#}\lambda ,\sigma _{D}\left( \lambda \right) \right) \right\} .
\label{setLambda}
\end{equation}%
We shall use the following easy fact:%
\begin{equation}
\Lambda +\left( \left\{ 0_{X^{\prime }}\right\} \times \mathbb{R}_{+}\right)
\text{ is a convex cone in }X^{\prime }\times \mathbb{R}\text{.}  \label{3}
\end{equation}

\begin{lemma}
\label{lemma1}If $D$ is a (nonempty convex) closed set and $B$ is nonempty,
then%
\begin{equation*}
w^{\ast }-\limfunc{cl}\left( \Lambda +\left\{ 0_{X^{\prime }}\right\} \times
\mathbb{R}_{+}\right) =\limfunc{epi}\sigma _{B}.
\end{equation*}
\end{lemma}

\textbf{Proof.} We first prove that
\begin{equation}
\Lambda \subset \limfunc{epi}\sigma _{B}.  \label{4}
\end{equation}

For each $\left( \lambda ,x\right) \in \beta \left( D\right) \times B$ one
has (since $\mathbb{A}x\in D$)%
\begin{equation*}
\left\langle \mathbb{A}^{\#}\lambda ,x\right\rangle =\left\langle \lambda ,%
\mathbb{A}x\right\rangle \leq \sigma _{D}\left( \lambda \right)
\end{equation*}%
and, consequently,
\begin{equation*}
\sigma _{B}\left( \mathbb{A}^{\#}\lambda \right) \leq \sigma _{D}\left(
\lambda \right) \text{ and }\left( \mathbb{A}^{\#}\lambda ,\sigma _{D}\left(
\lambda \right) \right) \in \limfunc{epi}\sigma _{B}.
\end{equation*}%
Then, (\ref{4}) holds and one has
\begin{equation*}
\Lambda +\left\{ 0_{X^{\prime }}\right\} \times \mathbb{R}_{+}\subset
\limfunc{epi}\sigma _{B}+\left\{ 0_{X^{\prime }}\right\} \times \mathbb{R}%
_{+}=\limfunc{epi}\sigma _{B}.
\end{equation*}%
Since $\limfunc{epi}\sigma _{B}$ is $w^{\ast }-$closed, we obtain the
inclusion $\left[ \subset \right] $ in Lemma \ref{lemma1}. We prove now the
reverse inclusion.

Let $\left( x^{\prime },s\right) \in \left( X^{\prime }\times \mathbb{R}%
\right) \diagdown w^{\ast }-\limfunc{cl}\left( \Lambda +\left\{ 0_{X^{\prime
}}\right\} \times \mathbb{R}_{+}\right) .$ By Hahn-Banach Theorem, there
exists $\left( a,r\right) \in X\times \mathbb{R}$ such that

\begin{equation}
\left\langle x^{\prime },a\right\rangle +rs>\sup\limits_{\left( \lambda
,\theta \right) \in \beta \left( D\right) \times \mathbb{R}_{+}}\left[
\left\langle \mathbb{A}^{\#}\lambda ,a\right\rangle +r\left( \sigma
_{D}\left( \lambda \right) +\theta \right) \right] .  \label{5}
\end{equation}%
Since $0_{Y^{\prime }}\in \beta \left( D\right) ,$ we have%
\begin{equation*}
\left\langle x^{\prime },a\right\rangle +rs>\sup\limits_{\theta \geq
0}r\theta
\end{equation*}%
and, by this, $r\leq 0.$

We now observe that $\left( \beta \left( D\right) \right) ^{-}=\limfunc{rec}%
\left( D\right) $ (see, e.g., \cite[Corollary 3.3.6]{CHL23}), where
\begin{equation*}
\limfunc{rec}\left( D\right) =\left\{ z\in Y:w+t\left( z-w\right) \in
D,\forall \left( w,t\right) \in D\times \mathbb{R}_{+}\right\}
\end{equation*}%
denotes the recession cone of the nonempty closed convex set $D\subset Y.$

So, for $r=0$ we get, by (\ref{5}),
\begin{equation*}
\left\langle x^{\prime },a\right\rangle >\sup\limits_{\lambda \in \beta
\left( D\right) }\left\langle \lambda ,\mathbb{A}a\right\rangle =\delta _{%
\limfunc{rec}\left( D\right) }\left( \mathbb{A}a\right) .
\end{equation*}%
We have $\mathbb{A}a\in \limfunc{rec}\left( D\right) $ and $\left\langle
x^{\prime },a\right\rangle >0.$ Picking $\overline{x}\in B,$ it follows that%
\begin{equation*}
\mathbb{A}\left( \overline{x}+ta\right) =\mathbb{A}\overline{x}+t\mathbb{A}%
a\in D,\forall t\geq 0,
\end{equation*}%
so that $\sigma _{B}\left( x^{\prime }\right) \geq \sup\limits_{t\geq
0}\left\langle x^{\prime },\overline{x}+ta\right\rangle =+\infty$, and $%
\left( x^{\prime },s\right) \notin \limfunc{epi}\sigma _{B}.$

Assume now that $r<0.$ By (\ref{5}) we have%
\begin{equation*}
\begin{array}{ll}
\left\langle x^{\prime },-\frac{a}{r}\right\rangle -s & >\sup\limits_{%
\lambda \in \beta \left( D\right) }\left( -\left\langle \mathbb{A}%
^{\#}\lambda ,\frac{a}{r}\right\rangle -\sigma _{D}\left( \lambda \right)
\right) \medskip \\
& =\sup\limits_{\lambda \in \beta \left( D\right) }\left( \left\langle
\lambda ,-\mathbb{A}\left( \frac{a}{r}\right) \right\rangle -\sigma
_{D}\left( \lambda \right) \right) \medskip =\delta _{D}^{\ast \ast }\left( -%
\mathbb{A}\left( \frac{a}{r}\right) \right) =\delta _{D}\left( -\mathbb{A}%
\left( \frac{a}{r}\right) \right) .%
\end{array}%
\end{equation*}%
Therefore, $\mathbb{A}\left( -\frac{a}{r}\right) \in D,$ $-\frac{a}{r}\in B,$
$\sigma _{B}\left( x^{\prime }\right) \geq \left\langle x^{\prime },-\frac{a%
}{r}\right\rangle >s,$ and, $\left( x^{\prime },s\right) \notin \limfunc{epi}%
\sigma _{B}.$\hfill $\square \medskip $

Let us introduce the set
\begin{equation}  \label{setKK}
K:=\Lambda +\limfunc{epi}\sigma _{C}\subset X^{\prime }\times \mathbb{R},
\end{equation}%
where $\Lambda =\bigcup\limits_{\lambda \in \beta \left( D\right) }\left\{
\left( \mathbb{A}^{\#}\lambda ,\sigma _{D}\left( \lambda \right) \right)
\right\} .$

The proof of the next lemma involving $K$ requires some more notation. Given
a nonempty subset $A$ of a linear space,\ we denote by $\limfunc{co}A$ and
by $\limfunc{cone}A$ the convex hull of $A$ and the convex cone containing
the origin generated by $A$, respectively. Obviously,%
\begin{equation*}
\limfunc{cone}A=\mathbb{R}_{+}\limfunc{co}A=\limfunc{co}\left( \mathbb{R}%
_{+}A\right) .
\end{equation*}

\begin{lemma}
\label{lemma2}Let $C$ be a nonempty convex set. Then,
\begin{equation*}
K=\limfunc{cone}\left( \limfunc{epi}\sigma _{C}\cup \Lambda \right) .
\end{equation*}
\end{lemma}

\textbf{Proof.} We first prove that $K$ is a convex cone. One has%
\begin{equation*}
K=\Lambda +\left[ \left( \left\{ 0_{X^{\prime }}\right\} \times \mathbb{R}%
_{+}\right) +\limfunc{epi}\sigma _{C}\right] =\left[ \Lambda +\left( \left\{
0_{X^{\prime }}\right\} \times \mathbb{R}_{+}\right) \right] +\limfunc{epi}%
\sigma _{C}.
\end{equation*}%
By (\ref{3}), $K$ is the Minkowski sum of two convex cones, that is, a
convex cone. Since $\left( 0_{X^{\prime }},0\right) \in \Lambda \cap
\limfunc{epi}\sigma _{C},$ the convex cone $K$ contains $\Lambda \cup
\limfunc{epi}\sigma _{C}$ and, by this,%
\begin{equation*}
\limfunc{cone}\left( \Lambda \cup \limfunc{epi}\sigma _{C}\right) \subset K.
\end{equation*}
Moreover,
\begin{equation*}
\Lambda +\limfunc{epi}\sigma _{C}\subset \limfunc{cone}\left( \Lambda \cup
\limfunc{epi}\sigma _{C}\right) ,
\end{equation*}%
and we are done.\hfill $\square \medskip $

We denote by $\Gamma \left( X\right) $ the set of lower semicontinuous
convex proper functions on $X.$ Recall that $B:=\mathbb{A}^{-1}\left(
D\right) .$\medskip

\begin{lemma}
\label{lemma3}Assume that $C$ and $D$ are (nonempty convex) closed sets and $%
B\cap C\neq \emptyset .$ Then we have%
\begin{equation*}
\limfunc{epi}\sigma _{B\cap C}=\overline{K}:=w^{\ast }-\limfunc{cl}K.
\end{equation*}
\end{lemma}

\textbf{Proof.} Note that $\delta _{B}$ and $\delta _{C}$ are $\Gamma \left(
X\right) -$functions such that $\func{dom}\delta _{B}\cap \func{dom}\delta
_{C}\neq \emptyset .$ Consequently, from \cite[(3)]{BJ05},%
\begin{equation*}
\limfunc{epi}\sigma _{B\cap C}=\limfunc{epi}\left( \delta _{B}+\delta
_{C}\right) ^{\ast }=w^{\ast }-\limfunc{cl}\left( \limfunc{epi}\sigma _{B}+%
\limfunc{epi}\sigma _{C}\right) .
\end{equation*}%
By Lemmas \ref{lemma1} and \ref{lemma2}, we have%
\begin{equation*}
\begin{array}{ll}
\limfunc{epi}\sigma _{B\cap C} & =w^{\ast }-\limfunc{cl}\left( \overline{%
\Lambda +\left\{ 0_{X^{\prime }}\right\} \times \mathbb{R}_{+}}+\limfunc{epi}%
\sigma _{C}\right) \medskip \\
& =w^{\ast }-\limfunc{cl}\left[ \left( \Lambda +\left\{ 0_{X^{\prime
}}\right\} \times \mathbb{R}_{+}\right) +\limfunc{epi}\sigma _{C}\right]
\medskip \\
& =w^{\ast }-\limfunc{cl}\left[ \Lambda +\left( \left\{ 0_{X^{\prime
}}\right\} \times \mathbb{R}_{+}+\limfunc{epi}\sigma _{C}\right) \right]
\medskip \\
& =w^{\ast }-\limfunc{cl}\left( \Lambda +\limfunc{epi}\sigma _{C}\right) =%
\overline{K},%
\end{array}%
\end{equation*}%
which completes the proof.\hfill $\square \medskip $

\begin{lemma}
\label{lemma4}Let $f\in \Gamma \left( X\right) $ and $C$ and $D$ be
(nonempty convex) closed sets such that $B\cap C\cap \func{dom}f\neq
\emptyset .$ Then,%
\begin{equation*}
\limfunc{epi}\left( f+\delta _{B\cap C}\right) ^{\ast }=w^{\ast }-\limfunc{cl%
}\left( \limfunc{epi}f^{\ast }+K\right) .
\end{equation*}
\end{lemma}

\textbf{Proof.} Note that $f$ and $\delta _{B\cap C}$ are $\Gamma \left(
X\right) -$functions such that $\func{dom}f\cap \func{dom}\delta _{B\cap
C}\neq \emptyset .$ It follows that%
\begin{equation*}
\limfunc{epi}\left( f+\delta _{B\cap C}\right) ^{\ast }=w^{\ast }-\limfunc{cl%
}\left( \limfunc{epi}f^{\ast }+\limfunc{epi}\sigma _{B\cap C}\right) .
\end{equation*}%
By Lemma \ref{lemma3},%
\begin{equation*}
\begin{array}{ll}
\limfunc{epi}\left( f+\delta _{B\cap C}\right) ^{\ast } & =w^{\ast }-%
\limfunc{cl}\left( \limfunc{epi}f^{\ast }+\overline{K}\right) \medskip \\
& =w^{\ast }-\limfunc{cl}\left( \limfunc{epi}f^{\ast }+K\right) ,%
\end{array}%
\end{equation*}%
and we are done. \hfill $\square \medskip $

Recall that \textrm{(A)} means that $f\left( x\right) \geq 0$ is consequence
of the inclusion $x\in C\cap \mathbb{A}^{-1}\left( D\right) =C\cap B.$%
\medskip

\begin{proposition}[Existence theorem]
\label{prop1}Assume that $C$ and $D$ are (nonempty convex) closed sets. The
following statements are equivalent:

\textrm{(i)} $C\cap B\neq \emptyset ,$

\textrm{(ii)} $\left( 0_{X^{\prime }},-1\right) \notin \overline{K}.$
\end{proposition}

\textbf{Proof.} $\left[ \mathrm{(i)}\Longrightarrow \mathrm{(ii)}\right] $
Since $\left( 0_{X^{\prime }},-1\right) \notin \limfunc{epi}\sigma _{B\cap
C},$ Lemma \ref{lemma3} yields $\left( 0_{X^{\prime }},-1\right) \notin
\overline{K}.$

$\left[ \mathrm{(ii)}\Longrightarrow \mathrm{(i)}\right] $ By Hahn-Banach
Theorem there exists $\left( a,r\right) \in X\times \mathbb{R}$ such that
\begin{equation}
\begin{array}{ll}
-r & =\left\langle 0_{X^{\prime }},a\right\rangle -r\medskip \\
& >\sup\limits_{\left( \lambda ,x^{\prime },\theta \right) \in \beta \left(
D\right) \times \beta \left( C\right) \times \mathbb{R}_{+}}\left[
\left\langle \mathbb{A}^{\#}\lambda +x^{\prime },a\right\rangle +r\left(
\sigma _{D}\left( \lambda \right) +\sigma _{C}\left( x^{\prime }\right)
+\theta \right) \right] .%
\end{array}
\label{6}
\end{equation}

Taking $\left( \lambda ,x^{\prime },0\right) =\left( 0_{Y^{\prime
}},0_{X^{\prime }},0\right) $ we obtain $-r>0$ and, dividing by $-r$ in (\ref%
{6}), we get%
\begin{equation*}
\begin{array}{ll}
1 & >\sup\limits_{\lambda \in \beta \left( D\right) }\left( \left\langle
\mathbb{A}^{\#}\lambda ,-\frac{a}{r}\right\rangle -\sigma _{D}\left( \lambda
\right) \right) +\sup\limits_{x^{\prime }\in \beta \left( C\right) }\left(
\left\langle x^{\prime },-\frac{a}{r}\right\rangle -\sigma _{C}\left(
x^{\prime }\right) \right) \medskip \\
& =\sigma _{D}^{\ast }\left( \mathbb{A}\left( -\frac{a}{r}\right) \right)
+\sigma _{C}^{\ast }\left( -\frac{a}{r}\right) \medskip \\
& =\delta _{D}\left( \mathbb{A}\left( -\frac{a}{r}\right) \right) +\delta
_{C}\left( -\frac{a}{r}\right) ,%
\end{array}%
\end{equation*}%
and, finally, $-\frac{a}{r}\in B\cap C,$ which is nonvoid. \hfill $\square
\medskip $

\begin{corollary}
\label{Cor3}Assume that $D$ is a (nonempty convex) closed set. Then $B$ is
nonempty if and only if%
\begin{equation*}
\left( 0_{X^{\prime }},-1\right) \notin w^{\ast }-\limfunc{cl}\left(
\bigcup\limits_{\lambda \in \beta \left( D\right) }\left\{ \left( \mathbb{A}%
^{\#}\lambda ,\sigma _{D}\left( \lambda \right) \right) \right\} +\left(
\left\{ 0_{X^{\prime }}\right\} \times \mathbb{R}_{+}\right) \right) .
\end{equation*}
\end{corollary}

\textbf{Proof.} Apply Proposition \ref{prop1} for $C=X.$ \hfill $\square
\medskip $

We now state the $w^{\ast }-$ closedness characterization for Farkas type
lemma.\medskip

\begin{theorem}[3rd characterization of Farkas' Lemma]
\label{Th3}Let $f\in \Gamma \left( X\right) ,$ $C$ (resp. $D$) be a
(nonempty convex) closed subset of $X$ (resp. $Y$) such that $C\cap \mathbb{A%
}^{-1}\left( D\right) \cap \func{dom}f\neq \emptyset .$ Then, the statements

$(\mathcal{A})$\ \ \ \ $x\in C\cap \mathbb{A}^{-1}(D)\ \ \Longrightarrow
f(x)\geq 0$

\noindent and

$(\mathcal{B})$ \ \ $\left\{
\begin{array}{ll}
& \exists \left( u^{\prime },v^{\prime },\lambda \right) \in \func{dom}%
f^{\ast }\times \beta \left( C\right) \times \beta \left( D\right) \ \text{
such\ that } \\
& f^{\ast }\left( u^{\prime }\right) +\sigma _{C}\left( v^{\prime }\right)
+\sigma _{D}\left( \lambda \right) \leq 0\text{ and }u^{\prime }+v^{\prime
}=-\mathbb{A}^{\#}\lambda%
\end{array}%
\right. $

\noindent are equivalent if and only if
\begin{equation*}
\limfunc{epi}f^{\ast }+K\text{ is }w^{\ast }-\text{closed regarding }\left\{
\left( 0_{X^{\prime }},0\right) \right\} .\mathrm{\ }
\end{equation*}
\end{theorem}

\textbf{Proof.} \textit{Necessity.} Assume that $(\mathcal{A})
\Longleftrightarrow (\mathcal{B})$ (in fact, $(\mathcal{A}) \Longrightarrow (%
\mathcal{B})$) and let $\left( 0_{X^{\prime }},0\right) \in \overline{%
\limfunc{epi}f^{\ast }+K}.$ By Lemma \ref{lemma4}, $\left( 0_{X^{\prime
}},0\right) \in \limfunc{epi}\left( f+\delta _{B\cap C}\right) ^{\ast },$
that means $(\mathcal{A})$ holds. By hypothesis, $(\mathcal{B})$ holds and
there exists $\left( u^{\prime },v^{\prime },\lambda \right) \in \func{dom}%
f^{\ast }\times \beta \left( C\right) \times \beta \left( D\right) $ such
that
\begin{equation*}
\left( 0_{X^{\prime }},0\right) =\left( u^{\prime },-\sigma _{C}\left(
v^{\prime }\right) -\sigma _{D}\left( \lambda \right) \right) +\left(
v^{\prime },\sigma _{C}\left( v^{\prime }\right) \right) +\left( \mathbb{A}%
^{\#}\lambda ,\sigma _{D}\left( \lambda \right) \right) .
\end{equation*}%
Thus,%
\begin{equation*}
\left( 0_{X^{\prime }},0\right) \in \limfunc{epi}f^{\ast }+\limfunc{epi}%
\sigma _{C}+\left\{ \left( \mathbb{A}^{\#}\lambda ,\sigma _{D}\left( \lambda
\right) \right) \right\} \subset \limfunc{epi}f^{\ast }+K.
\end{equation*}%
Hence, $\limfunc{epi}f^{\ast }+K$ is $w^{\ast }-$closed regarding $\left\{
\left( 0_{X^{\prime }},0\right) \right\} .$\textrm{\ }

{\textit{Sufficiency. }} Assume that $\limfunc{epi}f^{\ast }+K$ is $w^{\ast
}-$closed regarding $\left\{ \left( 0_{X^{\prime }},0\right) \right\} .$%
\textrm{\ }

One has to prove that $(\mathcal{A}) \Longleftrightarrow (\mathcal{B})$ or
just $(\mathcal{A}) \Longrightarrow (\mathcal{B})$ since $(\mathcal{B})
\Longleftrightarrow (\mathcal{A})$ always holds. By $(\mathcal{A})$ we have $%
\left( 0_{X^{\prime }},0\right) \in \limfunc{epi}\left( f+\delta _{B\cap
C}\right) ^{\ast }$ and, by Lemma \ref{lemma4}, $\left( 0_{X^{\prime
}},0\right) \in \overline{\limfunc{epi}f^{\ast }+K}.$ From the hypothesis we
obtain that $\left( 0_{X^{\prime }},0\right) \in \limfunc{epi}f^{\ast }+K$
and there exists $\left( u^{\prime },s\right) \in \limfunc{epi}f^{\ast },$ $%
\left( v^{\prime },t\right) \in \limfunc{epi}\sigma _{C}$ and $\lambda \in
\beta \left( D\right) $ such that%
\begin{equation*}
u^{\prime }+v^{\prime }+\mathbb{A}^{\#}\lambda =0_{X^{\prime }}\text{ and }%
s+t+\sigma _{D}\left( \lambda \right) =0.
\end{equation*}%
Finally,%
\begin{equation*}
f^{\ast }\left( u^{\prime }\right) +\sigma _{C}\left( v^{\prime }\right)
+\sigma _{D}\left( \lambda \right) \leq s+t+\sigma _{D}\left( \lambda
\right) =0
\end{equation*}%
and $(\mathcal{B})$ holds.\hfill $\square \medskip $

\begin{remark}
Consider the graph of $f^{\ast }\in \Gamma \left( X\right) ,$ namely%
\begin{equation*}
\func{gph}f^{\ast }:=\left\{ \left( x^{\prime },f^{\ast }\left( x^{\prime
}\right) \right) :x^{\prime }\in \func{dom}f^{\ast }\right\} .
\end{equation*}%
We have%
\begin{equation*}
\begin{array}{ll}
\limfunc{epi}f^{\ast }+K & =\limfunc{epi}f^{\ast }+\limfunc{epi}\sigma
_{C}+\Lambda \medskip \\
& =\left[ \func{gph}f^{\ast }+\left\{ 0_{X^{\prime }}\right\} \times \mathbb{%
R}_{+}\right] +\limfunc{epi}\sigma _{C}+\Lambda \medskip \\
& =\func{gph}f^{\ast }+\left[ \left\{ 0_{X^{\prime }}\right\} \times \mathbb{%
R}_{+}+\limfunc{epi}\sigma _{C}\right] +\Lambda \medskip \\
& =\func{gph}f^{\ast }+\limfunc{epi}\sigma _{C}+\Lambda \medskip \\
& =\func{gph}f^{\ast }+K.%
\end{array}%
\end{equation*}
\end{remark}

\begin{corollary}
\label{Cor4}Let $\emptyset \neq C\subset X,$ $\emptyset \neq D\subset Y,$
and $\left( x^{\prime },r\right) \in X^{\prime }\times \mathbb{R}.$ Consider
the following statements:

$(\mathcal{A}_1)$\ \ \ \ $x\in C\cap \mathbb{A}^{-1}\left( D\right)
\Longrightarrow \left\langle x^{\prime },x\right\rangle \leq r,$

$(\mathcal{B}_{3})$ \ \ $\exists \lambda \in \beta \left( D\right) $ such
that $\sigma _{C}\left( x^{\prime }-\mathbb{A}^{\#}\lambda \right) +\sigma
_{D}\left( \lambda \right) \leq r.$

\noindent We have $(\mathcal{B}_3) \Longrightarrow (\mathcal{A}_1)$. If,
moreover, $C$ and $D$ are closed and convex, and $C\cap \mathbb{A}%
^{-1}\left( D\right) \neq \emptyset ,$ then
\begin{equation*}
\left[(\mathcal{A}_1) \Longleftrightarrow (\mathcal{B}_3) \right]
\Longleftrightarrow K\text{ is }w^{\ast }-\text{ closed regarding }\left\{
\left( x^{\prime },r\right) \right\} .
\end{equation*}
\end{corollary}

\textbf{Proof. }Apply Theorem \ref{Th3} with $f\left( x\right)
=r-\left\langle x^{\prime },x\right\rangle ,$ noticing that
\begin{equation*}
\limfunc{epi}f^{\ast }+K=\func{gph}f^{\ast }+K=-\left\{ \left( x^{\prime
},r\right) \right\} +K
\end{equation*}%
and that being $w^{\ast }-$closed regarding $\left\{ \left( 0_{X^{\prime
}},0\right) \right\} $ amounts to saying that $K$ is $w^{\ast }-$closed
regarding $\left\{ \left( x^{\prime },r\right) \right\} .$\hfill \hfill $%
\square \medskip $

\begin{corollary}
\label{Cor5}Let $\left( x^{\prime },r\right) \in X^{\prime }\times \mathbb{R}
$, $C$ be a (nonempty convex) closed subset of $X$ such that $C\cap \mathbb{A%
}^{-1}\left( Q\right) \neq \emptyset $ and the convex cone $Q$ be closed in $%
Y.$ Then, the statements

$(\mathcal{A}_{2})$\ \ \ \ $x\in C\cap \mathbb{A}^{-1}\left( Q\right)
\Longrightarrow \left\langle x^{\prime },x\right\rangle \leq r$

\noindent and

$(\mathcal{B}_4) $ \ \ $\exists \lambda \in \beta \left( Q\right) =Q^{-} $
such that $\sigma _{C}\left( x^{\prime }-\mathbb{A}^{\#}\lambda \right) \leq
r$

\noindent are equivalent if and only if%
\begin{equation*}
K=\limfunc{epi}\sigma _{C}+\mathbb{A}^{\#}\left( Q^{-}\right) \times \left\{
0\right\} \mathrm{\ }\text{is }w^{\ast }-\text{closed regarding }\left\{
\left( x^{\prime },r\right) \right\} .
\end{equation*}
\end{corollary}

\textbf{Proof.} It is an immediate consequence of Corollary \ref{Cor4}%
.\hfill $\square \medskip $

\begin{corollary}
\label{Cor6}Let $a^{\prime }\in X^{\prime }$ and the (nonempty) convex cone $%
Q$ be closed in $Y.$ Then, the statements

$(\mathcal{A}_{3})$\ \ \ \ $\mathbb{A}x\in Q\Longrightarrow \left\langle
a^{\prime },x\right\rangle \geq 0$

\noindent and

$(\mathcal{B}_5)$\ \ \ \ $\exists \mu \in Q^{+}$ such that $\mathbb{A}%
^{\#}\mu =a^{\prime }$

\noindent are equivalent if and only if
\begin{equation*}
\mathbb{A}^{\#}\left( Q^{+}\right) \text{ is }w^{\ast }-\text{closed
regarding }\left\{ a^{\prime }\right\} .
\end{equation*}
\end{corollary}

\textbf{Proof.} Apply Corollary \ref{Cor5} with $C=X,$ $D=Q$ and $\left(
x^{\prime },r\right) =-\left( a^{\prime },0\right) .$ We have that the
following statements,\newline
\textrm{(i) }$\mathbb{A}x\in Q\Longrightarrow \left\langle a^{\prime
},x\right\rangle \geq 0$\newline
and\newline
\textrm{(ii) }$\exists \lambda \in \beta \left( Q\right) =Q^{-}$ such that $%
a^{\prime }=-\mathbb{A}^{\#}\lambda ,$ \newline
are equivalent if and only if\newline
\textrm{(iii) }$K=\left\{ 0_{X^{\prime }}\right\} \times \mathbb{R}_{+}+%
\mathbb{A}^{\#}\left( Q^{-}\right) \times \left\{ 0\right\} $ is $w^{\ast }-$%
closed regarding $\left\{ \left( -a^{\prime },0\right) \right\} .$

Noticing that \textrm{(ii)} reads $a^{\prime }\in \mathbb{A}^{\#}\left(
Q^{+}\right) $ while \textrm{(iii)} reads equivalently\newline
\textrm{(iii.a) }$\mathbb{A}^{\#}\left( Q^{-}\right) \times \mathbb{R}_{+}$
is $w^{\ast }-$closed regarding $\left\{ \left( -a^{\prime },0\right)
\right\} ,$\newline
\textrm{(iii.b) }$\mathbb{A}^{\#}\left( Q^{+}\right) \times \mathbb{R}_{+}$
is $w^{\ast }-$closed regarding $\left\{ \left( a^{\prime },0\right)
\right\} ,$\newline
and\newline
\textrm{(iii.c) }$\mathbb{A}^{\#}\left( Q^{+}\right) $ is $w^{\ast }-$closed
regarding $\left\{ a^{\prime }\right\}.$

\noindent We are done.\hfill $\square $\medskip

\begin{proposition}[Characterization of stable Farkas' lemma]
\label{prop2} Let $f\in \Gamma \left( X\right) ,$ $C$ (resp. $D$) be a
(nonempty convex) closed subset of $X$ (resp. $Y$) such that $C\cap \mathbb{A%
}^{-1}\left( D\right) \cap \func{dom}f\neq \emptyset .$ Then, the statements

$(\mathcal{A}_{4})$ $\forall \left( x^{\prime },\QTR{frametitle}{r}\right)
\in X^{\prime }\times \mathbb{R},$ $x\in C\cap \mathbb{A}^{-1}\left(
D\right) \Longrightarrow f\left( x\right) \geq \langle x^{\prime },x\rangle +%
\QTR{frametitle}{r}$

\noindent and

$(\mathcal{B}_6)$ $\left\{ \!\!\!\!\!\!
\begin{array}{ll}
& \forall \left( x^{\prime },\QTR{frametitle}{r}\right) \in X^{\prime
}\times \mathbb{R},\exists \left( u^{\prime },v^{\prime },\lambda \right)
\in \left( \func{dom}f^{\ast }-x^{\prime }\right) \times \beta \left(
C\right) \times \beta \left( D\right) \text{such that} \\
& u^{\prime}+v^{\prime }=-\mathbb{A}^{\#}\lambda \text{ and } f^{\ast
}\left( u^{\prime }+x^{\prime }\right) +\sigma _{C}\left( v^{\prime }\right)
+\sigma _{D}\left( \lambda \right) \leq -\QTR{frametitle}{r}%
\end{array}
\right. $

\noindent are equivalent if and only if
\begin{equation*}
\limfunc{epi}f^{\ast }+K\text{ is }w^{\ast }-\text{closed}.\mathrm{\ }
\end{equation*}
\end{proposition}

\textbf{Proof.} Observe that, given $\left( x^{\prime },\QTR{frametitle}{r}%
\right) \in X^{\prime }\times \mathbb{R},$ $\left( f-x^{\prime }-%
\QTR{frametitle}{r}\right) ^{\ast }\left( u^{\prime }\right) =%
\QTR{frametitle}{r}+f^{\ast }\left( u^{\prime }+x^{\prime }\right) $ for all
$u^{\prime }\in X^{\prime },$ $\func{dom}\left( f-x^{\prime }-%
\QTR{frametitle}{r}\right) ^{\ast }=\func{dom}f^{\ast }-x^{\prime },$ and $%
\limfunc{epi}\left( f-x^{\prime }-{r}\right) ^{\ast }=\limfunc{epi}f^{\ast
}-\left( x^{\prime },-{r}\right) .$ Hence, by Theorem \ref{Th3} applied to
the function $f-x^{\prime }- {r}$, the statements

$(\mathcal{A}_{4}^{rx^\prime })$ \ \ $x\in C\cap \mathbb{A}^{-1}\left(
D\right) \Longrightarrow f\left( x\right) \geq x^{\prime }+ r,$ and

$(\mathcal{B}_{6}^{rx^\prime })$ $\left\{
\begin{array}{ll}
& \exists \left( u^{\prime },v^{\prime },\lambda \right) \in \left( \func{dom%
}f^{\ast }-x^{\prime }\right) \times \beta \left( C\right) \times \beta
\left( D\right) \text{ such\ that } \\
& u^{\prime }+v^{\prime }=-\mathbb{A}^{\#}\lambda \text{ and } f^{\ast
}\left( u^{\prime }+x^{\prime }\right) +\sigma _{C}\left( v^{\prime }\right)
+\sigma _{D}\left( \lambda \right) \leq -\QTR{frametitle}{r},%
\end{array}
\right. $

\noindent are equivalent if and only if $\limfunc{epi}f^{\ast }+K$ is $%
w^{\ast }-$closed regarding $\left\{ \left( x^{\prime },-\QTR{frametitle}{r}%
\right) \right\} .$ Therefore, $(\mathcal{A}_{4})$ and $(\mathcal{B}_{6})$
are equivalent if and only if $\limfunc{epi}f^{\ast }+K$ is $w^{\ast }-$%
closed regarding $\left\{ \left( x^{\prime },-\QTR{frametitle}{r}\right)
\right\} $ for all $\left( x^{\prime },-\QTR{frametitle}{r}\right) \in
X^{\prime }\times \mathbb{R},$ that is, if and only if $\limfunc{epi}f^{\ast
}+K$ is $w^{\ast }-$closed.\hfill $\square \medskip $

Taking $X$ and $Y\ $normed linear spaces, $C=X,$ and $D$ a closed convex
cone in $Y,$ one gets the stable strong Farkas' lemma for linear operators
from $X$ to $Y$ \cite[Corollary 3.5]{JKL08} and \cite[Corollary 3.3]{JL08}.

\section{Characterization of a convex Farkas lemma for linear infinite
systems}

Given a possibly infinite family of continuous linear functionals $\left\{
a_{t}^{\prime }\right\} _{t\in T}\subset \left( X^{\prime }\right) ^{T},$ a
product of compact intervals $\prod\limits_{t\in T}\left[ \alpha _{t},\beta
_{t}\right] \subset \mathbb{R}^{T},$ a nonempty convex subset $C$ of $X,$
and a proper convex function $f\in \overline{\mathbb{R}}^{X},$ we are
concerned with the characterization of a Farkas type lemma involving the
consequence relation%
\begin{equation}
\left\{
\begin{array}{c}
\alpha _{t}\leq \left\langle a_{t}^{\prime },x\right\rangle \leq \beta
_{t},t\in T \\
x\in C%
\end{array}%
\right\} \Longrightarrow f\left( x\right) \geq 0.  \label{7}
\end{equation}

We are in the case that $Y=\mathbb{R}^{T}$ equipped with the product
topology and null vector $0_{T},$ $Y^{\prime }=\mathbb{R}^{\left( T\right) }$
(the subspace of $\mathbb{R}^{T}$ formed by the functions with a finite
support), $D=$ $\prod\limits_{t\in T}\left[ \alpha _{t},\beta _{t}\right] ,$
and $\mathbb{A}:X\longrightarrow \mathbb{R}^{T}$ such that $\mathbb{A}%
x=\left( \left\langle a_{t}^{\prime },x\right\rangle \right) _{t\in T}$ for
all $x\in X.$

\begin{remark}
\label{rem4.1} Given $\lambda =\left( \lambda _{t}\right) _{t\in T}\in
\mathbb{R}^{\left( T\right) },$\ we have%
\begin{equation*}
\begin{array}{ll}
\sigma _{D}\left( \lambda \right) & =\sup\limits_{\left( \eta _{t}\right)
_{t\in T}\in \tprod\limits_{t\in T}\left[ \alpha _{t},\beta _{t}\right]
}\sum\limits_{t\in T}\eta _{t}\lambda _{t}\medskip \\
& =\sum\limits_{t\in T,\alpha _{t}\leq s\leq \beta _{t}}\sup \lambda
_{t}\medskip s\ \ =\ \ \ \sum\limits_{t\in T}\left( \lambda _{t}^{+}\beta
_{t}-\lambda _{t}^{-}\alpha _{t}\right) ,%
\end{array}%
\end{equation*}%
where $\lambda _{t}^{+}=\max \left\{ \lambda _{t},0\right\} $ and $\lambda
_{t}^{-}=\max \left\{ -\lambda _{t},0\right\} .$ Observe that $\beta \left(
D\right) =\mathbb{R}^{\left( T\right) }$ due to the compactness of $%
\prod\limits_{t\in T}\left[ \alpha _{t},\beta _{t}\right] .$
\end{remark}

Since $\lambda _{t}\medskip =\lambda _{t}^{+}-\lambda _{t}^{-},$ we have $%
\mathbb{A}^{\#}\lambda =\sum\limits_{t\in T}\lambda _{t}a_{t}^{\prime
}=\sum\limits_{t\in T}\left( \lambda _{t}^{+}a_{t}^{\prime }-\lambda
_{t}^{-}a_{t}^{\prime }\right) $ and
\begin{equation*}
\Lambda =\bigcup\limits_{\lambda \in \mathbb{R}^{\left( T\right) }}\left\{
\sum\limits_{t\in T}\lambda _{t}^{+}\left( a_{t}^{\prime },\beta _{t}\right)
+\sum\limits_{t\in T}\lambda _{t}^{-}\left( -a_{t}^{\prime },-\alpha
_{t}\right) \right\} .
\end{equation*}

Let us consider the second moment cone $N$ (in the sense of \cite{DGL06}) of
the linear infinite system $\left\{ x\in X;\alpha _{t}\leq \left\langle
a_{t}^{\prime },x\right\rangle \leq \beta _{t},t\in T\right\} $ . We have%
\begin{equation}
\begin{array}{ll}
N\!\!\!\!\! & :=\limfunc{cone}\left[ \left\{ \left( a_{t}^{\prime },\beta
_{t}\right) ,t\in T\right\} \cup \left\{ \left( -a_{t}^{\prime },-\alpha
_{t}\right) ,t\in T\right\} \right] \medskip \\
& =\bigcup\limits_{\mu \in \mathbb{R}_{+}^{\left( T\right) }}\left\{
\sum\limits_{t\in T}\mu _{t}\left( a_{t}^{\prime },\beta _{t}\right)
\right\} +\bigcup\limits_{\eta \in \mathbb{R}_{+}^{\left( T\right) }}\left\{
\sum\limits_{t\in T}\eta _{t}\left( -a_{t}^{\prime },-\alpha _{t}\right)
\right\} .%
\end{array}
\label{11}
\end{equation}

\begin{proposition}[Relationship between $N$ and $\Lambda $]
\label{prop42} One has%
\begin{equation}
\Lambda \subset N\subset \Lambda +\left\{ 0_{X^{\prime }}\right\} \times
\mathbb{R}_{+}=N+\left\{ 0_{X^{\prime }}\right\} \times \mathbb{R}_{+}.
\label{10}
\end{equation}
\end{proposition}

\textbf{Proof.} Given $\lambda \in \mathbb{R}^{\left( T\right) },$ define $%
\mu =\left( \lambda _{t}^{+}\right) _{t\in T},\eta =\left( \lambda
_{t}^{-}\right) _{t\in T}\in \mathbb{R}_{+}^{\left( T\right) }.$ We obtain $%
\sum\limits_{t\in T}\mu _{t}\left( a_{t}^{\prime },\beta _{t}\right)
+\sum\limits_{t\in T}\eta _{t}^{-}\left( -a_{t}^{\prime },-\alpha
_{t}\right) \in N$ and, by this, $\Lambda \subset N.$

Given $\overline{t}\in T,$ define $\lambda \in \mathbb{R}^{\left( T\right) }$
such that $\lambda _{\overline{t}}=1$ and $\lambda _{t}=0$ otherwise. Then $%
\lambda \in \mathbb{R}_{+}^{\left( T\right) }$ and
\begin{equation*}
\left( a_{\overline{t}}^{\prime },\beta _{\overline{t}}\right)
=\sum\limits_{t\in T}\lambda _{t}^{+}\left( a_{t}^{\prime },\beta
_{t}\right) -\sum\limits_{t\in T}\lambda _{t}^{-}\left( a_{t}^{\prime
},\alpha _{t}\right) \in \Lambda .
\end{equation*}%
Similarly, setting $\lambda \in \mathbb{R}^{\left( T\right) }$ such that $%
\lambda _{\overline{t}}=-1$ and $\lambda _{t}=0$ otherwise, we obtain%
\begin{equation*}
-\left( a_{\overline{t}}^{\prime },\beta _{\overline{t}}\right)
=\sum\limits_{t\in T}\lambda _{t}^{+}\left( a_{t}^{\prime },\beta
_{t}\right) -\sum\limits_{t\in T}\lambda _{t}^{-}\left( a_{t}^{\prime
},\alpha _{t}\right) \in \Lambda .
\end{equation*}%
Therefore, $\left\{ \left( a_{t}^{\prime },\beta _{t}\right) ,t\in T\right\}
\cup \left\{ \left( -a_{t}^{\prime },-\alpha _{t}\right) ,t\in T\right\}
\subset \Lambda .$

Since $\Lambda +\left\{ 0_{X^{\prime }}\right\} \times \mathbb{R}_{+}$ is a
convex cone containing $\Lambda ,$ it follows that $N\subset \Lambda
+\left\{ 0_{X^{\prime }}\right\} \times \mathbb{R}_{+}.$ Hence,%
\begin{equation*}
N+\left\{ 0_{X^{\prime }}\right\} \times \mathbb{R}_{+}\subset \Lambda
+\left\{ 0_{X^{\prime }}\right\} \times \mathbb{R}_{+}\subset N+\left\{
0_{X^{\prime }}\right\} \times \mathbb{R}_{+},
\end{equation*}%
so that (\ref{10}) holds.\hfill $\square \medskip $

Now, statement $(\mathcal{B})$ reads
\begin{equation}
\begin{array}{c}
\exists \left( u^{\prime },v^{\prime },\lambda \right) \in \func{dom}f^{\ast
}\times \beta \left( C\right) \times \mathbb{R}^{\left( T\right) }\text{
such that }\medskip \\
u^{\prime }+v^{\prime }=-\sum\limits_{t\in T}\lambda _{t}a_{t}^{\prime }%
\text{ and }f^{\ast }\left( u^{\prime }\right) +\sigma _{C}\left( v^{\prime
}\right) \leq \sum\limits_{t\in T}\left( \lambda _{t}^{-}\alpha _{t}-\lambda
_{t}^{+}\beta _{t}\right) .%
\end{array}
\label{8}
\end{equation}

Applying Theorem \ref{Th1} we obtain\medskip

\begin{corollary}
\label{Cor7} The statements (\ref{7}) and (\ref{8}) are equivalent if and
only if%
\begin{equation*}
\mathbb{R}_{+}\left[ \bigcup\limits_{x\in C}\left\{ \left( x,\left(
\left\langle a_{t}^{\prime },x\right\rangle \right) _{t\in T},0\right)
\right\} +\!\!\!\bigcup\limits_{v\in \func{dom}f}\!\!\!\left\{ \left(
-v,0_{Y},f\left( v\right) \right) \right\} +\left\{ 0_{X}\right\} \times
\prod\limits_{t\in T}\left[ \alpha _{t},\beta _{t}\right] \times \mathbb{R}%
_{-}\right]
\end{equation*}%
is closed regarding $\left\{ \left( 0_{X},0_{T},-1\right) \right\} .$
\end{corollary}

Consider the statement%
\begin{equation}
\exists \lambda \in \mathbb{R}^{\left( T\right) }\text{ such that }f\left(
x\right) +\sum\limits_{t\in T}\lambda _{t}\left\langle a_{t}^{\prime
},x\right\rangle \geq \sum\limits_{t\in T}\left( \lambda _{t}^{+}\beta
_{t}-\lambda _{t}^{-}\alpha _{t}\right) ,\forall x\in X.  \label{9}
\end{equation}

Applying Theorem \ref{Th2} we obtain:

\begin{corollary}
\label{Cor8} Assume that $f \in \overline{\mathbb{R}}^{X}$ is convex proper
and $C\cap \func{dom}f\neq \emptyset .$ The statements (\ref{7}) and (\ref{9}%
) are equivalent if and only if%
\begin{equation*}
\mathbb{R}_{+}\left[ \bigcup\limits_{x\in C\cap \func{dom}f}\text{ }\left\{
\left( \left( \left\langle a_{t}^{\prime },x\right\rangle \right) _{t\in
T},f\left( x\right) \right) \right\} -\left( \prod\limits_{t\in T}\left[
\alpha _{t},\beta _{t}\right] \times \mathbb{R}_{-}\right) \right]
\end{equation*}
is closed regarding $\left\{ \left( 0_{T},-1\right) \right\} .$\medskip
\end{corollary}

Note that the convex cone $K$ defined in \eqref{setKK} reads here
\begin{equation*}
\widetilde{K}=\bigcup\limits_{\lambda \in \mathbb{R}^{\left( T\right)
}}\left\{ \left( \sum\limits_{t\in T}\lambda _{t}a_{t}^{\prime
},\sum\limits_{t\in T}\left( \lambda _{t}^{+}\beta _{t}-\lambda
_{t}^{-}\alpha _{t}\right) \right) \right\} +\func{epi}\sigma _{C}
\end{equation*}

Applying Theorem \ref{Th3} and Proposition \ref{prop2} we obtain
respectively, the following two consequences:

\medskip

\begin{corollary}
\label{Cor9}Assume that $f\in \Gamma \left( X\right) ,$ $C$ is a (nonempty
convex) closed set, and there exists $\overline{x}\in C\cap \func{dom}f$
such that $\alpha _{t}\leq \left\langle a_{t}^{\prime },\overline{x}%
\right\rangle \leq \beta _{t}$ for all $t\in T.$ Then, the statements (\ref%
{7}) and (\ref{8}) are equivalent if and only if
\begin{equation*}
\limfunc{epi}f^{\ast }+\widetilde{K}\text{ is }w^{\ast }-\text{closed
regarding }\left\{ \left( 0_{X^{\prime }},0\right) \right\} .
\end{equation*}
\end{corollary}

\begin{corollary}
\label{Cor10} Assume that $f\in \Gamma \left( X\right) ,$ $C$ is a (nonempty
convex) closed set, and there exists $\overline{x}\in C\cap \func{dom}f$
such that $\alpha _{t}\leq \left\langle a_{t}^{\prime },\overline{x}%
\right\rangle \leq \beta _{t}$ for all $t\in T.$ Then the statements

$(\mathcal{A}_{5})$\ \ \ \ $\forall \left( x^{\prime },\QTR{frametitle}{r}%
\right) \in X^{\prime }\times \mathbb{R},$ $x\in C,\alpha _{t}\leq
\left\langle a_{t}^{\prime },x\right\rangle \leq \beta _{t},t\in
T\Longrightarrow f\left( x\right) \geq \langle x^{\prime },x\rangle +%
\QTR{frametitle}{r}$

\noindent and

$(\mathcal{B}_{7})$ $\left\{ \!\!\!\!%
\begin{array}{ll}
& \exists \left( u^{\prime },v^{\prime },\lambda \right) \in \func{dom}%
(f^{\ast }-x^{\prime })\times \beta \left( C\right) \times \mathbb{R}%
^{\left( T\right) }\text{ such that } \\
& u^{\prime }+v^{\prime }=-\sum\limits_{t\in T}\lambda _{t}a_{t}^{\prime },%
\text{ }f^{\ast }\left( u^{\prime }+x^{\prime }\right) +\sigma _{C}\left(
v^{\prime }\right) \leq \sum\limits_{t\in T}\left( \lambda _{t}^{-}\alpha
_{t}-\lambda _{t}^{+}\beta _{t}\right) -r%
\end{array}%
\right. $

\noindent are equivalent if and only if
\begin{equation*}
\limfunc{epi}f^{\ast }+\widetilde{K}\text{ is }w^{\ast }-\text{closed.}
\end{equation*}
\end{corollary}

\section{Characterization of a concave Farkas lemma}

In this section, as in Section 3, $C$ (resp. $D$) is a nonempty closed
convex subset of $X$ (resp. $Y$), $\mathbb{A} : X \to Y$ a continuous linear
operator, and $f \in \Gamma (X)$. We are concerned with the statements

$(\mathcal{A}^\prime)$\ \ \ \ $x \in C \cap \mathbb{A}^{-1} (D) \ \
\Longrightarrow f(x) \leq 0$,

$(\mathcal{B}^{\prime })$ \ \ $\left\{
\begin{array}{ll}
& \forall x^{\prime }\in \limfunc{dom}f^{\ast },\exists (u^{\prime },\lambda
)\in \beta (C)\times \beta (D)\ \text{such\ that} \\
& x^{\prime }=u^{\prime }+\mathbb{A}^{\#}\lambda \text{ and }f^{\ast
}(x^{\prime })\geq \sigma _{C}(u^{\prime })+\sigma _{D}(\lambda ).%
\end{array}%
\right. $

\noindent One has straightforwardly $(\mathcal{B}^{\prime })\Longrightarrow (%
\mathcal{A}^{\prime })$ and we look for a characterization of the reverse
implication $(\mathcal{A}^{\prime })\Longrightarrow (\mathcal{B}^{\prime })$%
. To this end, as in Section 3, we consider
\begin{equation*}
\Lambda :=\bigcup\limits_{\lambda \in \beta \left( D\right) }\left( \mathbb{A%
}^{\#}\lambda ,\sigma _{D}\left( \lambda \right) \right) ,\ \ \ \ \ K:=%
\limfunc{epi}\sigma _{C}+\Lambda .
\end{equation*}%
Since $B:=\mathbb{A}^{-1}(D),$ we have by (\ref{4}) $\Lambda \subset
\limfunc{epi}\sigma _{B}$ and, by this,

\begin{equation*}
\Lambda \cup \limfunc{epi}\sigma _{C}\subset \limfunc{epi}\sigma _{B\cap C}.
\end{equation*}%
Since $\limfunc{epi}\sigma _{B\cap C}$ is a weak$^{\ast }$-closed convex
cone, it follows that $\text{w}^{\ast }-\text{cl cone}(\Lambda \cup \limfunc{%
epi}\sigma _{C})\subset \limfunc{epi}\sigma _{B\cap C}$, and, by Lemma \ref%
{lemma2}, that
\begin{equation}
\overline{K}:=w^{\ast }-\limfunc{cl}K\subset \limfunc{epi}\sigma _{B\cap C}.
\label{V1}
\end{equation}%
The statement $(\mathcal{A}^{\prime })$ reads $\delta _{B\cap C}\geq f$\ or,
equivalently,
\begin{equation}
\limfunc{epi}f^{\ast }\ \subset \ \limfunc{epi}\sigma _{B\cap C},  \label{V2}
\end{equation}%
while the statement $(\mathcal{B}^{\prime })$ reads straightforwardly
\begin{equation}
\limfunc{epi}f^{\ast }\ \subset \ K.  \label{V3}
\end{equation}%
Assuming $B\not=\emptyset $ we have (see Lemma \ref{lemma1})
\begin{equation*}
\limfunc{epi}\sigma _{B}=w^{\ast }-\limfunc{cl}\left( \Lambda
+\{0_{X^{\prime }}\}\times \mathbb{R}_{+}\right)
\end{equation*}%
and, provided $B\cap C\not=\emptyset $ (see Lemma \ref{lemma3}),
\begin{equation}
\limfunc{epi}\sigma _{B\cap C}=\overline{K}.  \label{V4}
\end{equation}

\begin{definition}
We will say that $K$ is \textit{{simili-$w^*$-closed} regarding $\limfunc{epi%
} f^*$ if
\begin{equation*}
\overline{K} \supset \limfunc{epi} f^* \ \Longrightarrow\ \ K \supset
\limfunc{epi} f^*.
\end{equation*}
}
\end{definition}

Note that
\begin{eqnarray*}
K\ \text{ is }w^{\ast }-\text{closed}\ \ &\Longrightarrow &K\ \ \text{is }%
w^{\ast }-\text{closed\ regarding}\ \limfunc{epi}f^{\ast } \\
&\Longrightarrow &K\ \text{is simili}-w^{\ast }-\text{closed\ regarding}\
\limfunc{epi}f^{\ast }.
\end{eqnarray*}

\begin{theorem}
\label{thmV1} Given $C \subset X$, $D \subset Y$ two nonempty closed convex
sets, $\mathbb{A}: X \to Y$ a continuous linear operator, and $f \in \Gamma
(X)$, consider the statements

$\mathrm{(i)}$ \ $(\mathcal{A}^{\prime })\Longleftrightarrow (\mathcal{B}%
^{\prime })$,

$\mathrm{(ii)}$ \ $K$ is ${simili}-w^*-{closed\ regarding} \limfunc{epi} f^*$%
.

\noindent We have $\mathrm{(i)} \Longrightarrow \mathrm{(ii)}$. If,
moreover, $C \cap \mathbb{A}^{-1} (D) \not= \emptyset$ then $\mathrm{(i)}
\Longleftrightarrow \mathrm{(ii)}$.
\end{theorem}

\textbf{Proof.} $\left[ \mathrm{(i)}\Rightarrow \mathrm{(ii)}\right] $ Let $%
\limfunc{epi}f^{\ast }\subset \overline{K}$. By \eqref{V1} we have $\limfunc{%
epi}f^{\ast }\subset \limfunc{epi}\sigma _{B\cap C}$ that is ($\mathcal{A}%
^{\prime }$) (see \eqref{V2}). By $(i)$ $(\mathcal{B}^{\prime })$ holds,
that is $\limfunc{epi}f^{\ast }\subset K$ (see \eqref{V3}).

$\left[ \mathrm{(ii)}\Rightarrow \mathrm{(i)}\right] $ Assume $B\cap
C\not=\emptyset $ and $(\mathcal{A}^{\prime })$ holds. By \eqref{V2} and %
\eqref{V4} we have $\limfunc{epi}f^{\ast }\subset \overline{K}$. By $\mathrm{%
(ii)}$ we obtain $\limfunc{epi}f^{\ast }\subset K,$ that is $(\mathcal{B}%
^{\prime })$ (see \eqref{V3}) and we are done. \hfill $\square $ \bigskip

\section{Applications}

\subsection{Application to constrained convex minimization problems}

With the same data as in Section 3, namely, $\mathbb{A}:X\rightarrow Y$ is a
linear continuous operator, $f\in \Gamma (X)$, $C\not=\emptyset $ ($%
D\not=\emptyset $ resp.) is a convex closed subset in $X$ ($Y$ resp.). For a
function $h\in \overline{\mathbb{R}}^{X}$, its subdifferential at $a\in X$
is
\begin{equation*}
\partial h(a):=\left\{
\begin{array}{ll}
\{x^{\ast }\in X^{\ast }\,:\,h(x)\geq h(a)+\langle x^{\ast },x-a\rangle
,\forall x\in X\}, & \text{if }h(a)\in \mathbb{R}, \\
\emptyset , & \text{else.}%
\end{array}%
\right.
\end{equation*}%
We denote by $N(C,a):=\partial \delta _{C}(a)$ the normal cone to $C$ at $%
a\in C.$

We now consider the problem
\begin{equation*}
\begin{array}{lll}
(P)\quad & \inf & f(x) \\
& \text{s.t.} & x \in C, \ \mathbb{A}(x) \in D.%
\end{array}%
\end{equation*}
Recall that
\begin{equation*}  \label{setK}
K:= \limfunc{epi}\sigma _{C} + \bigcup\limits_{\lambda \in \beta \left(
D\right) }\left\{ \left( \mathbb{A}^{\#}\lambda ,\sigma _{D}\left( \lambda
\right) \right) \right\}
\end{equation*}%
is a convex cone associated with the system of constraints $\{ \mathbb{A} x
\in D, x \in C\}$.

\begin{proposition}[Optimality condition for (P)]
\label{prop61nw} Assume that $\limfunc{epi}f^{\ast }+K$ is $w^{\ast}$-closed
regarding $\left\{ \left( 0_{X^{\prime }},0\right) \right\} $ and let $\bar
x \in C \cap \mathbb{A}^{-1} (D) \cap \func{dom} f $. Then the following
statements are equivalent:

\textrm{(i)} $\bar{x}$ is an optimal solution of \textrm{(P)},

\textrm{(ii)}$\ \ \left\{ \!\!\!\!\!\!\!%
\begin{array}{ll}
& \exists \left( u^{\prime },v^{\prime },\lambda \right) \in \func{dom}%
f^{\ast }\times \beta \left( C\right) \times \beta \left( D\right) \ \text{
such\ that} \\
& f^{\ast }\left( u^{\prime }\right) +\sigma _{C}\left( v^{\prime }\right)
+\sigma _{D}\left( \lambda \right) \leq -f(\bar{x})\text{ and }u^{\prime
}+v^{\prime }=-\mathbb{A}^{\#}\lambda ,%
\end{array}%
\right. $

\textrm{(iii)}\ \ \ \ $0_{X^{\prime }}\in \partial f(\bar{x})+N(C,\bar{x})+%
\mathbb{A}^{\#}\big(N(D,\mathbb{A}\bar{x})\big).$
\end{proposition}

\medskip

\textbf{Proof.} $[ \mathrm{(i)} \Longrightarrow \mathrm{(ii)}]$ Let $%
\alpha:= \inf \text{(P)} = f(\bar x) \in \mathbb{R}$. We have
\begin{equation}  \label{eq61}
x \in C \cap \mathbb{A}^{-1} (D)\ \ \ \Longrightarrow \ \ \ f(x)- \alpha
\geq 0,
\end{equation}
which means that $(\mathcal{A})$ in Theorem \ref{Th3} holds with $f$
replacing by $f - \alpha$. By this, \eqref{eq61} is equivalent to the
existence of $\left( u^{\prime },v^{\prime },\lambda \right) \in \func{dom}%
f^{\ast }\times \beta \left( C\right) \times \beta \left( D\right)$ such
that $u^{\prime }+v^{\prime }=-\mathbb{A}^{\#}\lambda$ and $(f -
\alpha)^{\ast }\left( u^{\prime }\right) +\sigma _{C}\left( v^{\prime
}\right) +\sigma _{D}\left( \lambda \right) \leq 0.$ This last statement
amounts to (ii).

$[ \mathrm{(ii)} \Longrightarrow \mathrm{(iii)}]$ By (ii) we have
\begin{equation*}
\big[ f ^{\ast }\left( u^{\prime }\right) + f(\bar x) - \langle u^\prime,
\bar x\rangle \Big] + \big[\sigma _{C}\left( v^{\prime }\right) - \langle
v^\prime , \bar x\rangle \Big] + \big[ \sigma _{D}\left( \lambda \right) -
\langle \lambda, \mathbb{A} \bar x\rangle \big] \leq 0.
\end{equation*}
Since the three brackets are nonnegative, they are all equal to $0$. That
means $u^\prime \in \partial f (\bar x) $, $v^\prime \in N(C, \bar x)$, $%
\lambda \in N(D, \mathbb{A} \bar x)$ and we have
\begin{equation*}
0_{X^\prime} = u^\prime + v^\prime + \mathbb{A}^{\#} \lambda \in \partial
f(\bar x) + N(C, \bar x) + \mathbb{A} ^{\#}\big(N(D, \mathbb{A} \bar x)\big),
\end{equation*}
that is (iii).

$[\mathrm{(iii)}\Longrightarrow \mathrm{(i)}]$ Assume that (iii) holds.
Since $\partial f(\bar{x})+N(C,\bar{x})\subset \partial (f+\delta _{C})(\bar{%
x})$ and $\mathbb{A}^{\#}\big(N(D,\mathbb{A}\bar{x})\big)\subset N(\mathbb{A}%
^{-1}(D),\bar{x})$ we have
\begin{equation*}
0_{X^{\prime }}\in \partial (f+\delta _{C})(\bar{x})+N(\mathbb{A}^{-1}(D),%
\bar{x})\subset \partial \big(f+\delta _{C}+\delta _{\mathbb{A}^{-1}(D)}\big)%
(\bar{x}).
\end{equation*}%
Therefore, $0_{X^{\prime }}\in \partial \big(f+\delta _{C}+\delta _{\mathbb{A%
}^{-1}(D)}\big)(\bar{x})$ that is (i). \hfill $\square $\bigskip

We now consider the dual problem of (P) corresponding to the perturbation
function $F$ defined by \eqref{eq22nw}. The so-called Lagrangian associated
to $F$ is (see \cite{Rock74}), for $(x,u^{\prime },\lambda )\in X\times
X^{\prime }\times Y^{\prime }$,
\begin{eqnarray*}
L(x,(u^{\prime },\lambda )) &=&\inf\limits_{(u,y)\in X\times Y}\big(%
F(x,(u,y))+(u^{\prime },u)+\langle \lambda ,y\rangle \big) \\
&=&\left\{
\begin{array}{lll}
& -\big(f^{\ast }(u^{\prime })-\langle u^{\prime }+\mathbb{A}^{\#}\lambda
,x\rangle +\sigma _{D}(\lambda )\big),\ \  & \mathrm{if}\ x\in C, \\
& +\infty ,\ \ \ \ \ \ \ \ \ \ \ \ \ \  & \mathrm{else}.%
\end{array}%
\right.
\end{eqnarray*}%
Then the Lagrangian dual of (P) reads
\begin{equation*}
\text{(P$^{\prime }$)}\quad \sup\limits_{(u^{\prime },\lambda )\in X^{\prime
}\times Y^{\prime }}\inf\limits_{x\in X}L(x,(u^{\prime },\lambda )).\hskip%
2.7cm
\end{equation*}%
%
%
%
%
%
%
%
%
%
%
%
We have $\inf \mathrm{(P)}\geq \sup \mathrm{(P}^{\prime })$ and, more
explicitly,
\begin{equation*}
\text{(P$^{\prime }$)}\quad \sup_{\substack{ (u^{\prime },v^{\prime
},\lambda )\in X^{\prime }\times X^{\prime }\times Y^{\prime }  \\ u^{\prime
}+v^{\prime }=-\mathbb{A}^{\#}\lambda }}-\Big(f^{\ast }(u^{\prime })+\sigma
_{C}(v^{\prime })+\sigma _{D}(\lambda )\Big).
\end{equation*}%
We write $\inf \mathrm{(P)}=\max \mathrm{(P}^{\prime })$ when there exists $%
(u^{\prime },v^{\prime },\lambda )\in X^{\prime }\times Y^{\prime }\times
Y^{\prime }$ such that $u^{\prime }+v^{\prime }=-\mathbb{A}^{\#}\lambda $
and $\inf \mathrm{(P)}=-\Big(f^{\ast }(u^{\prime })+\sigma _{C}(v^{\prime
})+\sigma _{D}(\lambda )\Big)$. That is the so-called strong duality between
(P) and (P$^{\prime }$). Note that strong duality holds in particular in the
case when $\inf \mathrm{(P)}=-\infty $.

\begin{proposition}[Strong duality for (P)]
\label{prop62nw} Consider the statements

\textrm{(i)} \ $\inf \mathrm{(P)} = \max (\mathrm{P}^\prime)$,

\textrm{(ii)} \ $\func{epi}f^{\ast }+K$ is $w^{\ast }$-closed regarding $%
\{0_{X^{\prime }}\}\times \mathbb{R}$.

We have \textrm{(i)} $\Longrightarrow $ \textrm{(ii)}. If, moreover, $\inf
\mathrm{(P)}\not=+\infty ,$ then \textrm{(ii)} $\Longrightarrow $ \textrm{%
(i) }.
\end{proposition}

\medskip

We need the following lemma for the proof of Propositions \ref{prop62nw} and %
\ref{prop64nw}.

\begin{lemma}
\label{lem63nw} Let $B = \mathbb{A}^{-1} (D)$. We have
\begin{equation*}
\overline{\func{epi} f^\ast + K} \ \subset \func{epi} \left( f +
\delta_{B\cap C}\right)^\ast.
\end{equation*}
\end{lemma}

\textbf{Proof of Lemma \ref{lem63nw}.} Let $(u^{\prime },s)\in \func{epi}%
f^{\ast }$, $(v^{\prime },t)\in \func{epi}\sigma _{C}$ and $\lambda \in
\beta (D)$. For each $x\in B\cap C$ we have
\begin{eqnarray*}
\langle u^{\prime }+v^{\prime }+\mathbb{A}^{\#}\lambda ,x\rangle -f(x)
&=&\langle u^{\prime },x\rangle -f(x)+\langle v^{\prime },x\rangle +\langle
\lambda ,\mathbb{A}x\rangle \\
&\leq &f^{\ast }(u^{\prime })+\sigma _{C}(v^{\prime })+\sigma _{D}(\lambda
)\leq s+t+\sigma _{D}(\lambda ).
\end{eqnarray*}%
Taking the supremum over $x\in B\cap C$ we obtain that
\begin{equation*}
\left( f+\delta _{B\cap C}\right) ^{\ast }(u^{\prime }+v^{\prime }+\mathbb{A}%
^{\#}\lambda )\leq s+t+\sigma _{D}(\lambda ).
\end{equation*}%
Therefore,
\begin{equation*}
\func{epi}f^{\ast }+K\subset \func{epi}\left( f+\delta _{B\cap C}\right)
^{\ast }
\end{equation*}%
and, since $\func{epi}\left( f+\delta _{B\cap C}\right) ^{\ast }$ is $%
w^{\ast }$-closed, the lemma is proved. \hfill $\square $

\medskip

\textbf{Proof of Proposition \ref{prop62nw}.} $[\mathrm{(i)}\Longrightarrow
\mathrm{(ii)}]$ Let $(0_{X^{\prime }},s)\in \overline{\func{epi}f^{\ast }+K}$%
. By Lemma \ref{lem63nw} we have $(f+\delta _{B\cap C})^{\ast }(0_{X^{\prime
}})\leq s,$ that is $-\inf \mathrm{(P)}\leq s$. By (i) there exists $%
(u^{\prime },v^{\prime },\lambda )\in X^{\prime }\times X^{\prime }\times
Y^{\prime }$ such that $u^{\prime }+v^{\prime }=-\mathbb{A}^{\#}\lambda $
and $f^{\ast }(u^{\prime })+\sigma _{C}(v^{\prime })+\sigma _{D}(\lambda
)\leq s.$ By this inequality there exists $(s_{1},s_{2})\in \mathbb{R}^{2}$
such that $(u^{\prime },s_{1})\in \func{epi}f^{\ast }$, $(v^{\prime
},s_{2})\in \func{epi}\sigma _{C}$ and $s_{1}+s_{2}=s-\sigma _{D}(\lambda )$%
. We then have
\begin{equation*}
(0_{X^{\prime }},s)=(u^{\prime },s_{1})+(v^{\prime },s_{2})+\big(\mathbb{A}%
^{\#}\lambda ,\sigma _{D}(\lambda )\big)\in \func{epi}f^{\ast }+K.
\end{equation*}

$[\mathrm{(ii)}\Longrightarrow \mathrm{(i)}]$ Since (i) holds if $\inf
\mathrm{(P)}=-\infty ,$ we can assume that $\alpha :=\inf \mathrm{(P)}\in
\mathbb{R}$. Then \eqref{eq61} holds and applying Theorem \ref{Th3} with $f$
replaced by $f-\alpha $ we obtain the existence of $(u^{\prime },v^{\prime
},\lambda )\in \func{dom}(f-\alpha )^{\ast }\times \beta (C)\times \beta (D)$
such that $u^{\prime }+v^{\prime }=-\mathbb{A}^{\#}\lambda $ and $(f-\alpha
)^{\ast }(u^{\prime })+\sigma _{C}(v^{\prime })+\sigma _{D}(\lambda )\leq 0$%
. Then we have
\begin{equation*}
\alpha \leq -\left( f^{\ast }(u^{\prime })+\sigma _{C}(v^{\prime })+\sigma
_{D}(\lambda )\right) \leq \sup \mathrm{(P}^{\prime })\leq \inf \mathrm{(P)}%
=\alpha
\end{equation*}%
and, finally, $\inf (\mathrm{P})=\max (\mathrm{P}^{\prime })$. \hfill $%
\square $

The notion of stable strong duality for a conic optimization problem was
introduced in \cite[page 2]{JSDL05} to indicate that the strong duality of
such a problem holds under arbitrary linear perturbations. This desirable
property was characterized through a closedness condition \cite[Theorem 3.1]%
{JSDL05}. We now obtain its corresponding version for $(\mathrm{P}).$

\medskip

\begin{proposition}[Stable strong duality for (P)]
\label{prop64nw} Consider the statements

\textrm{(i)} For all $x^{\prime }\in X^{\prime }$,
\begin{equation*}
\inf\limits_{x\in C\cap \mathbb{A}^{-1}(D)}\big[f(x)-\langle x^{\prime
},x\rangle \big]=\!\!\max_{\substack{ (u^{\prime },v^{\prime },\lambda )\in
X^{\prime }\times X^{\prime }\times Y^{\prime }  \\ u^{\prime }+v^{\prime }=-%
\mathbb{A}^{\#}\lambda }}\!\!\!-\Big[f^{\ast }(u^{\prime }+x^{\prime
})+\sigma _{C}(v^{\prime })+\sigma _{D}(\lambda )\Big],
\end{equation*}
\indent {\rm (ii)} \ $\func{epi}f^{\ast }+K$ is $w^{\ast }$-closed.

We have $\mathrm{(i)}\Longrightarrow \mathrm{(ii)}$. If, moreover, $\inf
\mathrm{(P)}\not=+\infty ,$ then $\mathrm{(ii)}\Longrightarrow \mathrm{(i)}$.
\end{proposition}

\medskip

\textbf{Proof.} $[\mathrm{(i)}\Longrightarrow \mathrm{(ii)}]$ Let $%
(x^{\prime },s)\in \func{epi}(f+\delta _{B\cap C})^{\ast }$. By (i) there
exists $(u^{\prime },v^{\prime },\lambda )\in X^{\prime }\times X^{\prime
}\times Y^{\prime }$ such that $f^{\ast }(u^{\prime }+x^{\prime })+\sigma
_{C}(v^{\prime })+\sigma _{D}(\lambda )\leq s$ and $u^{\prime }+v^{\prime }=-%
\mathbb{A}^{\#}\lambda $. Thus $f^{\ast }(u^{\prime }+x^{\prime })+\sigma
_{C}(v^{\prime })\leq s-\sigma _{D}(\lambda )$ and, by this, there exists $%
(s_{1}.s_{2})\in \mathbb{R}^{2}$ such that
\begin{equation*}
(u^{\prime }+x^{\prime },s_{1})\in \func{epi}f^{\ast },\ \ (v^{\prime
},s_{2})\in \func{epi}\sigma _{C},\ \ s_{1}+s_{2}=s-\sigma _{D}(\lambda ).
\end{equation*}%
Consequently,
\begin{equation*}
(x^{\prime },s)=(u^{\prime }+x^{\prime },s_{1})+(v^{\prime },s_{2})+(\mathbb{%
A}^{\#}\lambda ,\sigma _{D}(\lambda ))\in \func{epi}f^{\ast }+K,
\end{equation*}%
and by this, $\func{epi}(f+\delta _{B\cap C})^{\ast }\subset \func{epi}%
f^{\ast }+K$. By Lemma \ref{lem63nw} we obtain that $\func{epi}f^{\ast }+K=%
\func{epi}(f+\delta _{B\cap C})^{\ast }$ is weak$^{\ast }$-closed.

$[\mathrm{(ii)}\Longrightarrow \mathrm{(i)}]$ Given $x^{\prime }\in
X^{\prime }$ define the $\Gamma (X)$-function $f_{x^{\prime
}}(x):=f(x)-\langle x^{\prime },x\rangle $ and consider the problem
\begin{equation*}
\mathrm{(P}_{x^{\prime }})\ \ \ \inf\limits_{x\in B\cap C}f_{x^{\prime }}(x).%
\hskip5cm
\end{equation*}%
Let us apply the statement $\mathrm{(ii)}\Longrightarrow \mathrm{(i)}$ in
Proposition \ref{prop62nw} with $f$ replaced by $f_{x^{\prime }}$. To this
end notice that, since $\inf \mathrm{(P)}\not=+\infty ,$ we have $\inf (%
\mathrm{P}_{x^{\prime }})\not=+\infty $. Now $(f_{x^{\prime }})^{\ast
}(u^{\prime })=f^{\ast }(u^{\prime }+x^{\prime })$ and
\begin{equation*}
\func{epi}(f_{x^{\prime }})^{\ast }+K=\{(x^{\prime },0)\}+\func{epi}f^{\ast
}+K\ \ \text{is }w^{\ast }-\mathrm{closed.}
\end{equation*}%
We then have, by Proposition \ref{prop62nw},%
\begin{equation*}
\begin{array}{ll}
\inf\limits_{x\in B\cap C}\big[f(x)-\langle x^{\prime },x\rangle \big]%
\!\!\!\! & =\inf \mathrm{(P}_{x^{\prime }}) \\
& =\!\max\limits_{\substack{ (u^{\prime },v^{\prime },\lambda )\in X^{\prime
}\times X^{\prime }\times Y^{\prime }  \\ u^{\prime }+v^{\prime }=-\mathbb{A}%
^{\#}\lambda }}\!\!\! -\Big[(f_{x^{\prime }})^{\ast }(u^{\prime })+\sigma
_{C}(v^{\prime })+\sigma _{D}(\lambda )\Big] \\
& =\!\max\limits_{\substack{ (u^{\prime },v^{\prime },\lambda )\in X^{\prime
}\times X^{\prime }\times Y^{\prime }  \\ u^{\prime }+v^{\prime }=-\mathbb{A}%
^{\#}\lambda }}\!\!\! -\Big[f^{\ast }(u^{\prime }+x^{\prime })+\sigma
_{C}(v^{\prime })+\sigma _{D}(\lambda )\Big]%
\end{array}%
\end{equation*}%
and we are done.\hfill $\square $

\subsection{Functional approximation\ by polynomials}

Consider the biobjective best approximation to $g\in \mathcal{C}\left( \left[
0,1\right] ,\mathbb{R}\right) $ by means of real polynomials of degree less
than $n$ from above for the $L_{\infty }$ (also called uniform) and for the $%
L_{1}$ norms on $\mathcal{C}\left( \left[ 0,1\right] ,\mathbb{R}\right) .$
Denoting by $x\in \mathbb{R}^{n}$ the decision variable, with associated
polynomial $p(x,t):=\sum_{i=1}^{n}x_{i}t^{i-1},$ the biobjective problem
reads%
\begin{equation*}
\begin{array}{lll}
\mathrm{(P)} \quad & \text{\textquotedblleft }\inf\nolimits_{x\in \mathbb{R}%
^{n}}\text{\textquotedblright } & \left( \max\nolimits_{t\in \left[ 0,1%
\right] }\left( \sum_{i=1}^{n}x_{i}t^{i-1}-g(t)\right) ,\sum_{i=1}^{n}\frac{%
x_{i}}{i}\right) \medskip \\
& \text{s.t.} & g(t)\leq \sum_{i=1}^{n}x_{i}t^{i-1},\ t\in \left[ 0,1\right].%
\end{array}%
\end{equation*}%
One of the ways to generate efficient solutions of $\mathrm{(P)}$ (an
approach called scalarization via $\varepsilon -$constraints in \cite[\S 4.1]%
{Ehrgott05}) consists of solving the linear semi-infinite parametric problem%
\begin{equation*}
\begin{array}{lll}
(\mathrm{P}_{\varepsilon })\quad & \inf\nolimits_{x\in \mathbb{R}^{n}} &
f(x):=\sum_{i=1}^{n}\frac{x_{i}}{i}\medskip \\
& \text{s.t.} & g(t)\leq \sum_{i=1}^{n}x_{i}t^{i-1}\leq g(t)+\varepsilon ,\
t\in \left[ 0,1\right] ,%
\end{array}%
\end{equation*}%
for different values of a parameter $\varepsilon >0.$ In the notation of
Section 4, $X=X^{\prime }=\mathbb{R}^{n}$ with null vector $0_{n},\ T=\left[
0,1\right] ,\ a_{t}^{\prime }=\left( 1,t,...,t^{n-1}\right) $, $\alpha
_{t}=g(t),$ and $\beta _{t}=g(t)+\varepsilon $ for all $t\in \lbrack 0,1].$
Moreover, $C=\mathbb{R}^{n}$ and $\limfunc{epi}\sigma _{C}=\left\{
0_{n}\right\} \times \mathbb{R}_{+}.$ By Remark \ref{rem4.1}, the set $%
\Lambda $ defined in \eqref{setLambda} becomes here
\begin{equation*}
\Lambda =\left\{ \left( \sum\limits_{t\in \left[ 0,1\right] }\lambda
_{t}\left( 1,t,...,t^{n-1}\right) ,\sum\limits_{t\in \left[ 0,1\right]
}\left( \lambda _{t}g(t)+\varepsilon \lambda _{t}^{+}\right) \right)
:\lambda \in \mathbb{R}^{\left( \left[ 0,1\right] \right) }\right\} .
\end{equation*}%
Denote by $N$ the second moment cone of the linear semi-infinite constraint
system of \textrm{(P}$_{\varepsilon }$) (recall (\ref{11})). Proposition \ref%
{prop42} gives that
\begin{equation*}
\Lambda +\left\{ 0_{n}\right\} \times \mathbb{R}_{+}=N+\left\{ 0_{n}\right\}
\times \mathbb{R}_{+}.
\end{equation*}%
By Proposition \ref{prop1}, $(\mathrm{P}_{\varepsilon })$ is consistent if
and only if $\left( 0_{n},-1\right) \notin \overline{\Lambda +\left\{
0_{n}\right\} \times \mathbb{R}_{+}}$ or, equivalently, $\left(
0_{n},-1\right) \notin \overline{N}$ (see \ \cite[Theorem 4.4]{GL98}).

Farkas' lemma for a feasible problem $(\mathrm{P}_{\varepsilon })$
characterizes the consequence relation\newline
\begin{equation}
\left\{ g(t)\leq \sum_{i=1}^{n}x_{i}t^{i-1}\leq g(t)+\varepsilon ,\ t\in
\left[ 0,1\right] \right\} \Longrightarrow \sum_{i=1}^{n}\frac{x_{i}}{i}\geq
\alpha ,  \label{eq62nww}
\end{equation}%
with $\alpha \in \mathbb{R},$ by the condition that (see Corollary \ref{Cor8}
with $f(x):=\sum_{i=1}^{n}\frac{x_{i}}{i}-\alpha $): there exists $\lambda
\in \mathbb{R}^{([0,1])}$ such that
\begin{equation*}
\sum_{i=1}^{n}x_{i}\left( \frac{1}{i}+\sum\limits_{t\in \left[ 0,1\right]
}\lambda _{t}t^{i-1}\right) -\sum\limits_{t\in \left[ 0,1\right] }\left(
\lambda _{t}g(t)+\varepsilon \lambda _{t}^{+}\right) \geq \alpha ,\forall
x\in \mathbb{R}^{n},
\end{equation*}%
that is to say, there exists $\lambda \in \mathbb{R}^{([0,1])}$ such that

\begin{equation}  \label{eq63nww}
\left\{
\begin{array}{ll}
& \sum\limits_{t\in \left[ 0,1\right]} \lambda_t t^{i-1} = - \frac{1}{i}, \
i = 1, \cdots, n, \\
& - \sum\limits_{t\in \left[ 0,1\right] } \left( \lambda
_{t}g(t)+\varepsilon \lambda _{t}^{+}\right) \geq \alpha .%
\end{array}
\right.
\end{equation}
By Corollary \ref{Cor8}, we have that $\eqref{eq62nww} \Longrightarrow %
\eqref{eq63nww} $ if only if
\begin{equation*}
\mathbb{R}_{+}\left[ \bigcup_{x \in \mathbb{R}^n} \Big\{ \Big(%
\left(p(x,t)\right)_{t \in [0,1]} , \sum_{i=1}^{n}\frac{x_{i}}{i}-\alpha %
\Big) \Big\} - \left( \prod\limits_{t\in \left[ 0,1\right] }\left[
g(t),g(t)+\varepsilon \right] \times \mathbb{R}_{-}\right) \right]
\end{equation*}%
is closed regarding $\left\{ \left(0_{\left[ 0,1\right] },-1\right) \right\}
.$

Applying Corollary \ref{Cor9} with $f(x):=\sum_{i=1}^{n}\frac{x_{i}}{i}%
-\alpha $ we obtain that $\eqref{eq62nww}\Longrightarrow \eqref{eq63nww}$ if
only if
\begin{equation*}
\{0_{n}\}\times \mathbb{R}_{+}+\bigcup_{\lambda \in \mathbb{R}^{([0,1])}}%
\Big\{\sum_{t\in \lbrack 0,1]}\left( \lambda _{t}(1,t,\cdots
,t^{n-1}),\lambda _{t}g(t)+\varepsilon \lambda _{t}^{+}\right) \Big\}
\end{equation*}%
is closed regarding $\left\{ -(1,2,...,\frac{1}{n},\alpha )\right\} .$

\medskip

\textbf{Funding.} This research was partially supported by Grant
PID2022-136399NB-C21 funded by MICIU/AEI/10.13039/501100011033 and by
ERDF/EU. \bigskip

\textbf{Data availability declaration.} The manuscript has no associated
data. \bigskip

\textbf{Declarations}\bigskip

\textbf{Competing Interest declaration.} There are no competing interests
related to this research.

\end{document}